\newif\ifspringer
\newif\ifelsevier
\springertrue 

\ifspringer
\RequirePackage{fix-cm} 
\documentclass[smallcondensed,nospthms]{svjour3}
\usepackage[paperwidth=15.51cm,paperheight=23.51cm,left=1.75cm,right=1.75cm,top=2cm,bottom=2cm]{geometry}
\fi

\ifelsevier
\documentclass[preprint,3p,times]{elsarticle}
\fi

\usepackage[english]{babel}
\usepackage[utf8]{inputenc}
\usepackage{pdflscape} 
\usepackage{amssymb}
\usepackage{amsmath}
\usepackage{amsthm}
\usepackage{mathtools}
\ifspringer
\usepackage{bm} 
\fi
\usepackage[normalem]{ulem}
\usepackage{enumitem}
\usepackage{marvosym}
\usepackage[dvipsnames,svgnames]{xcolor}
\usepackage{array}
\usepackage{multirow}
\ifspringer
\usepackage{mathptmx} 
\fi
\ifelsevier
\usepackage[nodots]{numcompress} 
\fi

\usepackage{graphics}
\usepackage{graphicx}
\ifspringer
\usepackage[FIGBOTCAP,TABTOPCAP,scriptsize]{subfigure} 
\fi
\ifelsevier
\usepackage[FIGBOTCAP,TABTOPCAP]{subfigure} 
\fi
\usepackage{tikz}
\usetikzlibrary{calc,shapes,arrows,backgrounds}
\tikzstyle{every picture}+=[remember picture]
\usepackage{pgfplots} 

\usepackage{epstopdf} 
\usepackage{hyperref}
\hypersetup{colorlinks=true,linkcolor=blue,citecolor=blue,filecolor=blue,urlcolor=blue}
\usepackage[all]{hypcap}

\mathtoolsset{showonlyrefs} 
\allowdisplaybreaks 
\ifspringer
\smartqed 
\fi
\ifelsevier
\biboptions{sort&compress,numbers}
\fi

\usepackage{xspace}
\makeatletter
\DeclareRobustCommand\onedot{\futurelet\@let@token\@onedot}
\newcommand{\@onedot}{\ifx\@let@token.\else.\null\fi\xspace}
 
\newcommand{\ie}{{i.e}\onedot}

\makeatother

\newtheorem{remark}{Remark}[section]
\newtheorem{definition}{Definition}[section]
\newtheorem{theorem}{Theorem}[section]
\newtheorem{lemma}{Lemma}[section]
\newtheorem{proposition}{Proposition}[section]
\newtheorem{corollary}{Corollary}[section]

\usepackage[mathscr]{euscript}
\usepackage{amsopn}

\newcommand\NN{\mathbb{{N}}}

\newcommand\RR{\mathbb{{R}}}
\newcommand\CC{\mathbb{{C}}}
\newcommand\ZZ{\mathbb{{Z}}}

\DeclareMathAlphabet{\mathpzc}{OT1}{pzc}{m}{it}

\newcommand\rc{{\rm c}}

\newcommand\rf{{\rm f}}

\newcommand\ri{{\rm i}}

\newcommand\e{{\rm e}}
\newcommand\Circ{{\rm circ}}

\newcommand\bb{{\bf b}}
\newcommand\bc{{\bf c}}

\newcommand\bd{{\bf d}}
\newcommand\be{{\boldsymbol e}}
\newcommand\bbf{{\bf f}}

\newcommand\bn{{\bf n}}

\newcommand\bq{{\bf q}}
\newcommand\br{{\bf r}}

\newcommand\bu{{\bf u}}

\newcommand\bx{{\boldsymbol x}}
\newcommand\by{{\bf y}}
\newcommand\bz{{\boldsymbol z}}

\newcommand\bD{{\bf D}}

\newcommand\bX{{\boldsymbol X}}

\newcommand\balpha{{\boldsymbol\alpha}}
\newcommand\bbeta{{\boldsymbol\beta}}
\newcommand\bdelta{{\boldsymbol \delta}}

\newcommand\bgamma{{\boldsymbol\gamma}}

\newcommand\bPhi{{\boldsymbol\Phi}}
\newcommand\bPsi{{\boldsymbol\Psi}}

\newcommand\cM{{\mathcal M}}

\newcommand\cE{{\mathcal E}}

\newcommand\cJ{{\mathcal J}}

\ifspringer
\journalname{\dots}
\fi

\makeatletter
\providecommand{\doi}[1]{%
  \begingroup
    \let\bibinfo\@secondoftwo
    \urlstyle{rm}%
    \href{http://dx.doi.org/#1}{%
      doi:\discretionary{}{}{}%
      \nolinkurl{#1}%
    }%
  \endgroup
}
\makeatother

\begin{document}
\newcommand{\titletext}{
Convergence and normal continuity analysis of non-stationary subdivision schemes near extraordinary vertices and faces
}
\newcommand{\titlerunningtext}{
Convergence and normal continuity of non-stationary subdivision
}

\newcommand{\abstracttext}{
Convergence and normal continuity analysis of a bivariate non-stationary (level-dependent) subdivision scheme for 2-manifold meshes with arbitrary topology is still an open issue. Exploiting ideas from the theory of asymptotically equivalent subdivision schemes, in this paper we derive new sufficient conditions for establishing convergence
and normal continuity of any rotationally symmetric, non-stationary, subdivision scheme near an extraordinary vertex/face.
}

\ifspringer
\newcommand{\separ}{\and}
\fi

\newcommand{\keytext}{
Non-stationary subdivision \separ  Extraordinary vertex/face \separ Convergence \separ Normal continuity
}

\newcommand{\MSCtext}{
26A15 \separ 68U07
}

\ifspringer
\title{\titletext
}
\titlerunning{\titlerunningtext} 

\author{
Costanza Conti \and
Marco Donatelli \and
Lucia Romani \and
Paola Novara
}
\authorrunning{C.~Conti, M.~Donatelli, L.~Romani, P.~Novara} 

\institute{
C.~Conti \at
Dipartimento di Ingegneria Industriale, Universit\`{a} di Firenze,
Viale Morgagni 40/44, 50134 Firenze, Italy\\
\email{costanza.conti@unifi.it}
\and
M.~Donatelli \at
Dipartimento di Scienza e Alta Tecnologia, Universit\`{a} dell'Insubria,
Via Valleggio 11, 22100 Como, Italy\\
\email{marco.donatelli@uninsubria.it}
\and
L.~Romani  (\Letter) \at
Dipartimento di Matematica, Alma Mater Studiorum Universit\`{a} di Bologna,
Piazza di Porta San Donato 5, 40126 Bologna, Italy\\
\email{lucia.romani@unibo.it}
\and
P.~Novara \at
Dipartimento di Scienza e Alta Tecnologia, Universit\`{a} dell'Insubria,
Via Valleggio 11, 22100 Como, Italy\\
\email{paola.novara@uninsubria.it}
}

\date{Received: date / Accepted: date}

\maketitle

\begin{abstract}
\abstracttext
\subclass{\MSCtext}
\end{abstract}
\fi


\section{Introduction}
This paper provides a general procedure to check convergence of non-stationary (level-dependent) subdivision schemes in the neighborhood of an extraordinary vertex/face. It also gives sufficient conditions for the limit surface to be normal continuous at the limit point of an extraordinary vertex/face.
To the best of our knowledge, the only contributions in this domain are the works in \cite{CP13,DLL92,JSD02}, where specific schemes are considered.
The difficulties concerning the analysis of a level-dependent subdivision scheme
in the neighborhood of an extraordinary vertex/face, are due to the fact that the well-established approach based on the spectral analysis of the subdivision matrix and on the study of the characteristic map is not applicable.
Thus, we use and generalize the notion of asymptotical equivalence between stationary and non-stationary subdivision schemes (known only for schemes defined on regular meshes), and show that normal continuity of a non-stationary scheme in the vicinity of an extraordinary element can be obtained by assuming
that the matrix sequence identifying it converges towards the matrix $S$ (identifying
a $C^1$-regular, standard, stationary scheme) faster than $\lambda_1^k$, where $\lambda_1$ denotes the real, double subdominant eigenvalue of $S$.
The sufficient conditions we propose are used for the analysis of the family of approximating non-stationary subdivision schemes presented in \cite{FMW14}. The members of the latter family are a generalization of exponential spline surfaces to quadrilateral meshes of arbitrary topology whose normal continuity is conjectured and shown only by numerical evidence in \cite[Section 5]{FMW14}.\\
Due to the lack of existing theoretical results for the analysis of level-dependent subdivision schemes, we believe that our
contribution could mark a first step forward towards a deeper understanding of non-stationary subdivision with a consequent increase of its  use in different fields of application.

\subsection{Motivation}
Non-stationary subdivision schemes were introduced more than 20 years ago with the aim of enriching the class of limit functions of stationary schemes and have very different and distinguished properties. Indeed, it is well-known that stationary subdivision schemes are not capable of generating circles, ellipses, or to deal with level-dependent tension parameters that allow the user to arbitrarily modify the shape of a subdivision limit.
Non-stationary schemes generate function spaces that are much richer. For example, in the univariate case, they include exponential B-splines or $C^\infty$ limits with bounded support as the Rvachev-type function (see, e.g., \cite{DL02}).
The generation capabilities of level-dependent schemes  (especially the capability of generating exponential-polynomials) is important in several applications, e.g., in biological imaging \cite{Badoual,ContiRomaniUnser2015,DTS12,DU13,uhlmann14}, in geometric design-approximation \cite{DTU12,JL91,LY10,RMS16,WW02} and in isogeometric analysis \cite{Hughes}. Moreover, level-dependent subdivision schemes include Hermite schemes that do not only model curves and surfaces, but also their gradient fields (such schemes are again considered of interest both in geometric modelling and biological imaging, see, e.g.,   \cite{conti16,conti17,ContiRomaniUnser2015,JY17,R10}).
Additionally, non-stationary subdivision schemes are at the base of non-stationary wavelet and frame constructions that, being level adapted, are certainly more flexible \cite{cotronei17,DHLR14,H12,VBU07}.
Unfortunately, in practice, the use of subdivision is mostly restricted to the class of stationary subdivision schemes even though
the non-stationary ones are equally simple to implement and highly intuitive in use: from an implementation point of view changing coefficients with the levels is not a crucial matter also in consideration of the fact that, in practice, only few subdivision iterations are performed.
On the contrary, a crucial limitation to the spread of level-dependent schemes, is a lack of general analysis methods, especially of methods for establishing their convergence and normal continuity. This motivates our study.

\subsection{Subdivision framework}
\begin{figure}[!h]
	\centering
	\resizebox{4cm}{!}{	\includegraphics{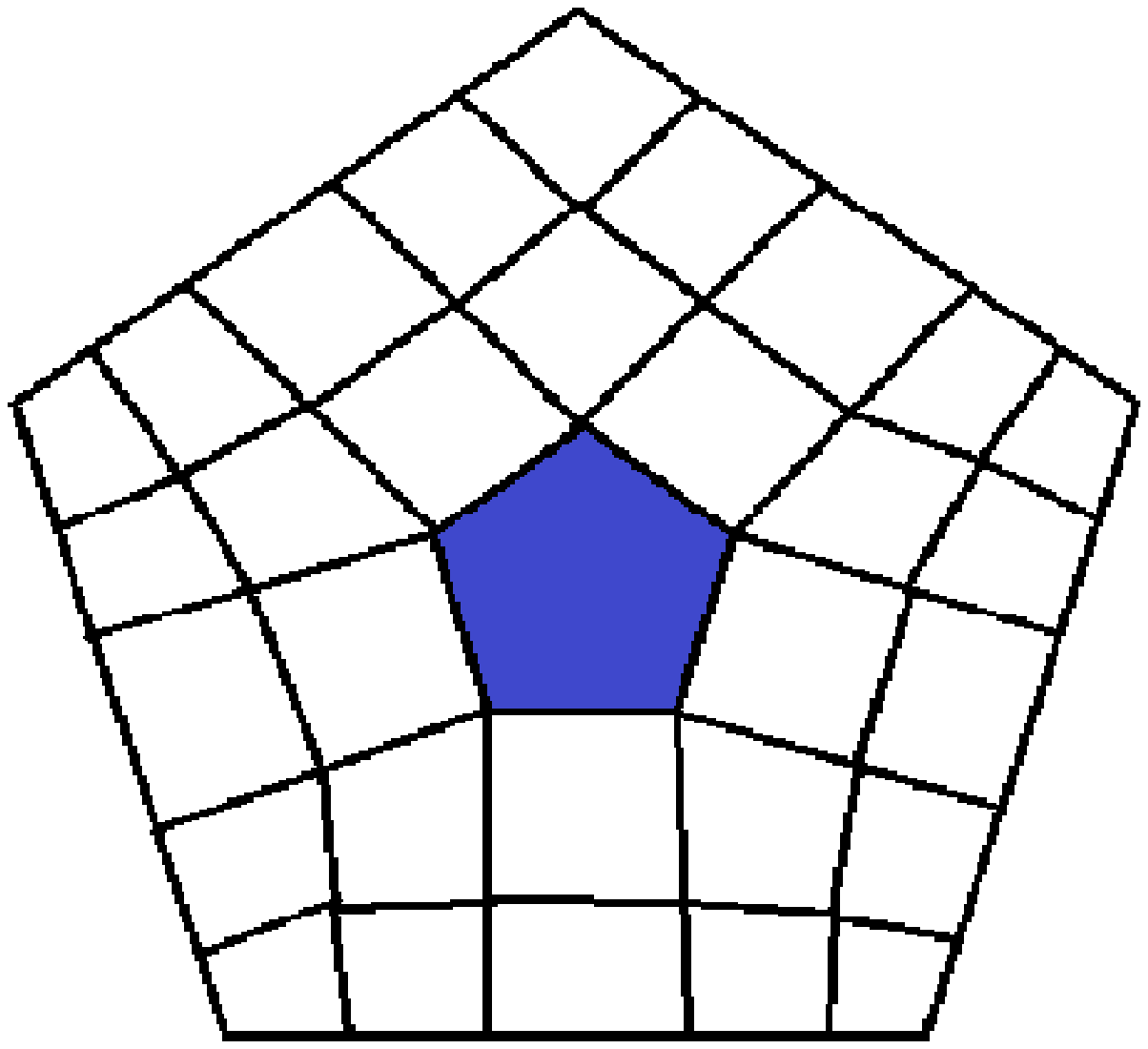}}
\resizebox{3.7cm}{!}{\includegraphics{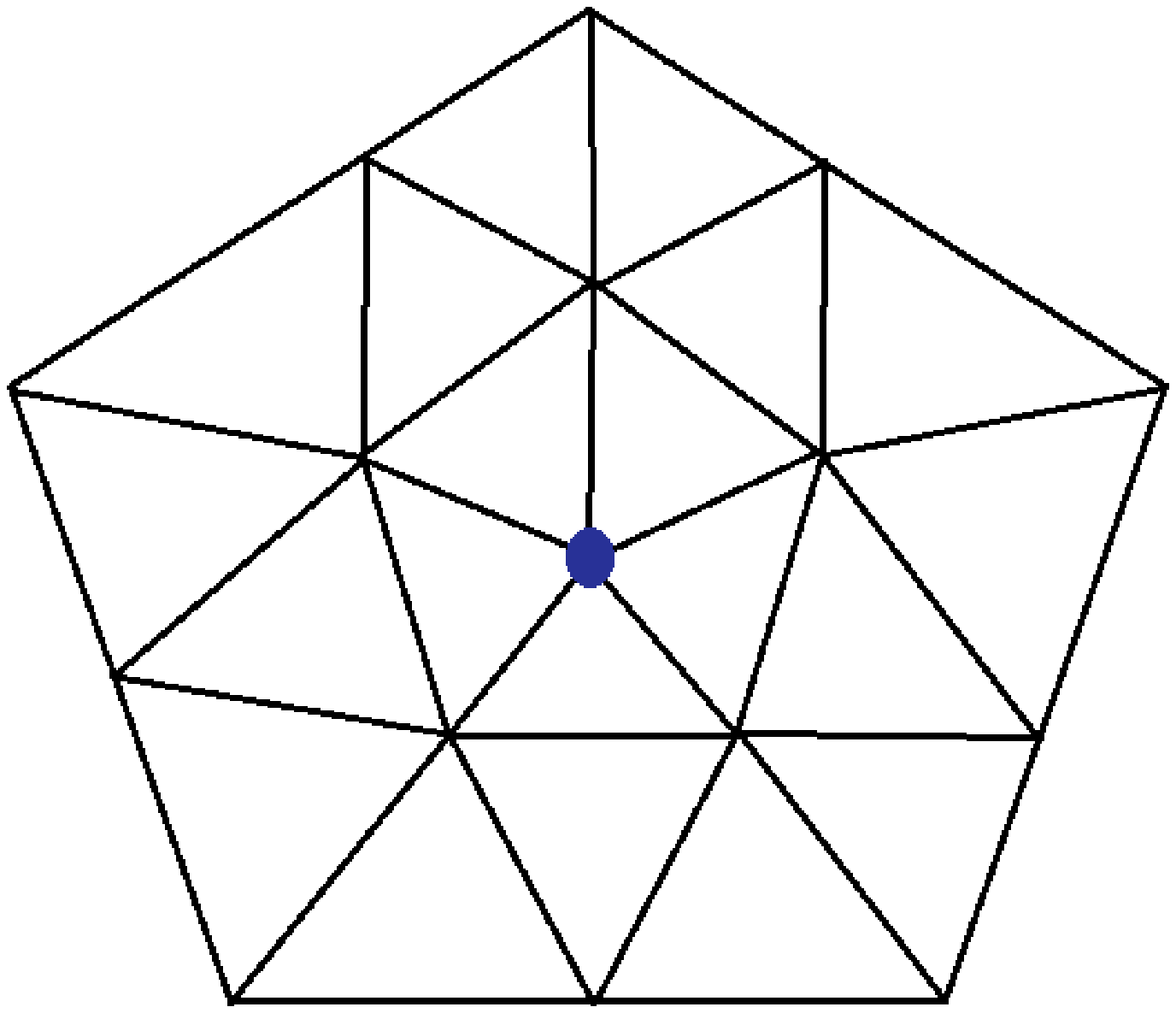}}
	\caption{Example of quadrilateral mesh containing an extraordinary face (left) and of triangular mesh containing an extraordinary vertex (right). (Color figure online.)}
\label{fig:extra_vert_face}
\end{figure}

Subdivision schemes are efficient iterative algorithms to produce smooth surfaces as the limit of a recursive process starting from a given coarse 2-manifold polygon mesh. A polygon mesh is considered to be 2-manifold if all its edges and faces are bounded, edges only intersect at vertices and are shared by at most two faces (boundary edges have one incident face, whereas inner edges have two incident faces); moreover, each vertex has either one connected
ring of faces around it (if inner) or one connected half-ring of faces (if boundary), see e.g. \cite{Gotsman}.
Each step of the recursive process produces a finer 2-manifold polygon mesh than the original one, containing many more vertices and polygonal faces.
The insertion of new vertices into a mesh requires modifications to both the topology (i.e., connectivity) and geometry (i.e., vertex
positions) of the mesh. For this reason, each subdivision scheme requires the specification of two rules: (i) a topologic refinement rule that describes how the connectivity of the mesh is to be modified in order to incorporate the new vertices being added to the mesh; and
(ii) a geometric refinement rule that describes how the geometry of the mesh is to be changed in order to accommodate
the new vertices being added (where these modifications may affect the position of previously-existing vertices).
The topologic and geometric refinement rules of a subdivision scheme may change with the refinement level or not. In the latter case the subdivision scheme is called \emph{stationary}, \emph{non-stationary}, or \emph{level-dependent}, otherwise.
Moreover, if the same set of geometric rules is used to determine all of the vertices within a single level of subdivision, the scheme is said to be \emph{uniform}.\\
Vertices and faces of a polygon mesh are classified by the so-called  \emph{vertex valence} and \emph{face valence}, respectively.
While the valence of a vertex is the number of edges incident to it, the valence of a face counts the number of edges that delimit it.
For a quadrilateral mesh, vertices and faces of valence 4 are called \emph{regular}. Differently,
for a triangular mesh regular vertices are the ones with valence 6, while regular faces have valence 3.
A \emph{regular mesh} or a \emph{regular region of a mesh} is a mesh/region where all vertices and faces are regular. Non-regular vertices/faces are called  \emph{extraordinary} (see Figure \ref{fig:extra_vert_face} for a graphical illustration of these two cases) and, whenever they appear, the mesh is said to be \emph{irregular} or of \emph{arbitrary topology}.
Accordingly, an \emph{irregular region of a mesh} contains \emph{extraordinary} vertices and/or faces. \\
A known analysis tool to investigate convergence and regularity of stationary subdivision schemes for regular meshes is the one based on symbols, originally proposed in \cite{CDM91,MP87} and successively exploited in \cite{DL02}.
To study convergence and regularity of a non-stationary subdivision scheme for regular meshes, Dyn and Levin \cite{DL95} proposed a method based on its comparison with a stationary scheme whose convergence and regularity are known.

\noindent In the case of meshes with arbitrary topology, we are currently able to study only convergence and regularity of stationary subdivision schemes near extraordinary vertices/faces, thanks to the results based on the spectral analysis of the subdivision matrix and on the study of the characteristic map \cite{PR08,R95,U2000,Z00,ZSS96}. However, in literature we can find no general results to analyze level-dependent  subdivision schemes near extraordinary elements.
To the best of our knowledge, the only contributions in this domain are the works in \cite{CP13,DLL92,JSD02}, where specific schemes are considered. Therefore, the  goal of our paper is to propose a general procedure to check if a non-stationary subdivision scheme is convergent in the neighborhood of an extraordinary vertex/face.  Moreover, it also aims at giving sufficient conditions for the limit surface to be normal continuous (in the sense of \cite{R95}) at the limit point of a sequence of extraordinary vertices/faces.

\medskip \noindent
The paper is organized as follows. In Section \ref{sec:2} we provide preliminaries  on bivariate, rotationally symmetric subdivision schemes. Then, in Section \ref{sec:new3}, we prove new results dealing with the $C^1$-convergence analysis of non-stationary subdivision schemes in regular regions. Next, in Section \ref{sec:3} (specifically, Subsection \ref{subsec:convergence}) sufficient conditions for proving convergence of a rotationally symmetric, non-stationary subdivision scheme near extraordinary vertices/faces are given. Finally, in Subsection \ref{sec:4} we also  give sufficient conditions to verify if the limit surface generated by a rotationally symmetric, convergent, non-stationary subdivision scheme is normal continuous at the extraordinary point, i.e., the limit of a sequence of extraordinary elements. Some application examples of the derived conditions are shown in Section \ref{sec:5}.

\section{Preliminaries on bivariate, rotationally symmetric subdivision schemes}\label{sec:2}

A bivariate subdivision scheme ${\mathscr S}$ is an iterative method that uses an initial polygonal mesh ${\mathcal M}^{(1)}$
to produce a sequence of denser and denser meshes $\{{\mathcal M}^{(k+1)}, k \in \NN\}$ that, when $k$ tends to infinity, converges to a smooth surface $\br$. In the sequel we use $k\geq 1$ instead of $k \in \NN$, omitting the trivial information that the refinement level is always assumed to be an integer.\\
Unless explicitly specified, we consider rotationally symmetric, local, uniform and non-stationary (level-dependent) subdivision schemes for meshes of arbitrary  topology, \ie,  subdivision schemes with symmetric, local refinement rules depending only on the level and eventually on the type of vertex (face point, edge point, vertex point in case of primal subdivision), but not on the vertex location.
We consider schemes that near extraordinary vertices/faces use rules that preserve their number and their location during refinements. This means that the number of extraordinary elements in ${\mathcal M}^{(1)}, \ldots,\ {\mathcal M}^{(k-1)}$, ${\mathcal M}^{(k)}$, ${\mathcal M}^{(k+1)}, \ldots$ remains unchanged.
The action of ${\mathscr S}$ in the regular regions of  $\cM^{(k)}$ can be described by the componentwise application of the refinement rules
\begin{equation} \label{eq:ss}
\rf_{\balpha}^{(k+1)}= \sum_{\bbeta \in \ZZ^2} \rc^{(k)}_{\balpha-2\bbeta}\ \rf^{(k)}_{\bbeta},\quad k\geq 1,  \;
\balpha \in \ZZ^2,
\end{equation}
where the set of coefficients $\bc^{(k)}=\{\rc^{(k)}_{\balpha},\ \alpha\in \ZZ^2\}$, also called the $k$-th level \emph{subdivision mask}, is finite due to the locality of the subdivision scheme. To simplify the analysis we also assume that all sets of coefficients have the same bounded supports.
Equivalently, the action of ${\mathscr S}$ on regular points  can be described by the application of
the subdivision operator ${\cal S}_{\bc^{(k)}}$, mapping componentwise the vector $\bbf^{(k)}$ into the corresponding vector of level $k+1$, i.e.,
$$
\bbf^{(k+1)}={\cal S}_{\bc^{(k)}} \bbf^{(k)}.
$$
The coefficients in \eqref{eq:ss} can be conveniently incorporated in the $k$-th level \emph{subdivision symbol}
$$
c^{(k)}(\bz)=\sum_{\balpha \in \ZZ^2}\rc^{(k)}_{\balpha} \bz^{\balpha}, \quad \bz \in (\CC \backslash \{0\})^2.
$$
The notation $\|{\cal S}_{\bc^{(k)}}\|_{\infty}$ is for the norm of the operator ${\cal S}_{\bc^{(k)}}$, i.e.,
\begin{equation}\label{normSO}
\|{\cal S}_{\bc^{(k)}}\|_{\infty}:=\max \left \{\sum_{\bbeta \in \ZZ^2} | \rc^{(k)}_{\balpha -2\bbeta}| \ : \ \balpha \in \{(0,0), (0,1), (1,0), (1,1) \} \right \}.
\end{equation}
In conclusion, when applied to regular regions, a subdivision scheme ${\mathscr S}$ can be equivalently identified with the sequence of subdivision operators $\{{\cal S}_{\bc^{(k)}}, \, k \geq 1 \}$, with the sequence of subdivision masks
$\{\bc^{(k)}, \, k \geq 1 \}$ or with the sequence of associated subdivision symbols $\{c^{(k)}(\bz), \, k \geq 1 \}$. \\
Instead, when applied to an irregular region, the subdivision rules relating the vertices of the $k$-th level mesh with those of the next level  $k+1$ are encoded in the rows of a non-singular local subdivision matrix $S_k$.
Thus, in the neighborhood of an extraordinary element the action of the subdivision scheme
${\mathscr S}$ is described by a sequence of non-singular local subdivision matrices $\{S_k, k \geq 1 \}$.

\begin{remark}
Note that the local subdivision matrix $S_k$ is also an alternative way to represent a subdivision step in regular regions.
\end{remark}

\smallskip \noindent
In the stationary setting we will use the notation $\bar{\mathscr S}$ to refer to a subdivision scheme that is not level-dependent. Hence, it will be identified with
\begin{itemize}
\item a subdivision operator, say ${\cal S}_{\bc}$, a subdivision mask
$\bc$ or an associated subdivision symbol $c(\bz)$, when applied to regular regions,
\item a local subdivision matrix $S$, when applied to an irregular region.
\end{itemize}

\subsection{Preliminaries for studying convergence of non-stationary subdivision schemes in regular regions}\label{sub:cr}

In the following, after recalling some well-known definitions, we present several useful results dealing with the convergence of a non-stationary subdivision scheme in regular regions (see, e.g., \cite{CRY16,DL95} for further details).

\begin{definition}
A non-stationary subdivision scheme ${\mathscr S}$ is called \emph{convergent} if, for any initial data $\bbf^{(1)} \in \ell(\ZZ^2)$, there exists a function $g_{\bbf^{(1)}} \in C(\RR^2)$ such that
$$
\lim_{\ell \rightarrow +\infty} \sup_{\balpha \in \ZZ^2} |g_{\bbf^{(1)}}( 2^{-\ell}\balpha)-{\rm f}^{(\ell+1)}_{\balpha}|=0,
$$
and if $g_{\bbf^{(1)}}$ is nonzero for at least one initial nonzero sequence $\bbf^{(1)}$.
For $r\ge 1 $, the subdivision scheme ${\mathscr S}$ is called \emph{$C^{r}$-convergent} if $g_{\bbf^{(1)}} \in C^{r}(\RR^2)$.
\end{definition}

\begin{definition}
For a convergent, stationary subdivision scheme $\bar{\mathscr S}:=\{{\cal S}_{\bc}\}$ the limit function obtained from the initial sequence $\bdelta=\{\delta_{{\bf 0},\balpha}, \, \balpha \in \ZZ^2\}$, denoted as
\begin{equation}\label{def:BLFstat}
\bar{\phi}:=\lim_{\ell \rightarrow +\infty} ({\cal S}_{\bc})^{\ell} \, \bdelta,
\end{equation}
is called the \emph{basic limit function} of the subdivision scheme.
\end{definition}

\begin{definition}
For a convergent, non-stationary subdivision scheme ${\mathscr S}:=\{{\cal S}_{\bc^{(\ell)}}, \, \ell \geq 1 \}$ the limit function
obtained from the initial sequence $\bdelta=\{\delta_{{\bf 0},\balpha}, \, \balpha \in \ZZ^2\}$, denoted as
\begin{equation}\label{def:BLFnonstat}
\phi_k:=\lim_{\ell \rightarrow +\infty} {\cal S}_{\bc^{(k+\ell)}} \, {\cal S}_{\bc^{(k+\ell-1)}} \, \ldots {\cal S}_{\bc^{(k)}} \,  \bdelta,
\end{equation}
is called the $k$-th member of the family of \emph{basic limit functions} $\{\phi_k, \, k\geq 1\}$ of the subdivision scheme.
\end{definition}

We remark that, in this paper, we consider non-stationary subdivision schemes that are \emph{non-singular} in the sense that they generate a zero limit if and only if the starting sequence is the zero sequence. Under this assumption we are guaranteed that, for each level $k \geq 1$, the shifts of the basic limit function $\phi_k$ are linearly independent \cite[Propostion 1.3]{CC13}.
The non-singularity assumption of the local subdivision matrix $S_k$ provides an analogous property in the neighborhood of an extraordinary element, in the sense that $S_k \bd_k={\bf 0}$ if and only if the control point vector is $\bd_k={\bf 0}$.

\begin{definition}\label{def:asympt_eq}
Let ${\mathscr S}$ and $\bar{\mathscr S}$ be subdivision schemes acting on regular regions with the subdivision masks $\{\bc^{(k)}\in \ell(\ZZ^2), \, k \geq 1\}$ and $\bc \in \ell(\ZZ^2)$, respectively.
If
\begin{equation}\label{eq:asympt_eq}
\sum_{k=1}^{+\infty}\|{\cal S}_{\bc^{(k)}}-{\cal S}_{\bc}\|_{\infty} < +\infty,
\end{equation}
then ${\mathscr S}$ and $\bar{\mathscr S}$ are said to be \emph{asymptotically equivalent} schemes.
\end{definition}

\begin{remark}\label{rem:CDMM}
As observed in \cite[page 2]{CDMM15}, \eqref{eq:asympt_eq} holds if and only if $$\displaystyle \sum_{k=1}^{+\infty}\|\bc^{(k)}-\bc\|_{\infty} < +\infty\quad \hbox{where} \quad \|\bc\|_{\infty}=\sup_{\balpha \in \ZZ^2} |\rc_{\balpha}|.$$
\end{remark}

\begin{theorem}\cite[Theorems 7-8 and Lemma 15]{DL95}\label{prop:asympt_equiv}
Let ${\mathscr S}$ and $\bar{\mathscr S}$ be asymptotically equivalent subdivision schemes acting on regular regions with the subdivision masks $\{\bc^{(k)} \in \ell(\ZZ^2), \, k \geq 1\}$ and $\bc \in \ell(\ZZ^2)$, respectively.
If $\bar{\mathscr S}$ is convergent, then ${\mathscr S}$ is also convergent and
$$
\lim_{k\rightarrow + \infty} \sup_{(u,v) \in \RR^2} | \phi_k(u,v) - \overline{\phi}(u,v)|=0,
$$
where $\overline{\phi}$ is the basic limit function of $\bar{\mathscr S}$ defined in \eqref{def:BLFstat} and
$\{\phi_k, k \geq 1\}$ the family of basic limit functions of ${\mathscr S}$ defined in \eqref{def:BLFnonstat}.
\end{theorem}

\begin{remark}
The condition of asymptotical equivalence in \eqref{eq:asympt_eq}, that guarantees convergence, could be relaxed by considering the fulfillment of the weaker condition of asymptotical similarity together with approximate sum rules of order 1, as shown in \cite{CCGP15}.
\end{remark}

\subsection{Preliminaries for studying convergence of rotationally symmetric, non-stationary  subdivision schemes in irregular regions}

We start our discussion by observing that we can restrict our analysis to a mesh $\cM^{(1)}$ with a single extraordinary element surrounded by a number of ``rings'' of ordinary vertices constituting the sub-mesh, here denoted by $\cE^{(1)}$, which determines a neighborhood of the extraordinary point, i.e., the limit of a sequence of extraordinary elements.
The number of rings and, consequently, the number of vertices in $\cE^{(1)}$ depends on the specific subdivision scheme (for example, there are $3$ ``rings''  in case of Loop's scheme). Obviously, the regular part of $\cM^{(1)}$ will be simply given by $\cM^{(1)}\setminus \cE^{(1)}$.

\smallskip In the neighborhood of an extraordinary vertex/face, each step of a subdivision algorithm can be conveniently encoded in the rows of a local subdivision matrix ${S}_k$ relating the vertices of the $k$-th level mesh with those of the next level.
The matrix ${S}_k$ has a different structure depending on the scheme properties and on the kind of extraordinary element (face or vertex) appearing in the $k$-th level mesh. Precisely, if the scheme is rotationally symmetric and the mesh contains an extraordinary face of valence $n$, in view of the fact that the valence-$n$ extraordinary face is surrounded by $n$ sectors, each composed by $p$ vertices, the local subdivision matrix ${S}_k$ is of the form
\begin{equation}\label{eq:Sk_block_dual}
{S}_{k}=\left(\begin{array}{cccc}
{B}_{0,k} & {B}_{1,k}  & \cdots & {B}_{n-1,k}\\
{B}_{n-1,k} & {B}_{0,k} & \cdots & {B}_{n-2,k}\\
\vdots & \vdots  & \ddots & \vdots\\
{B}_{1,k} & \cdots & {B}_{n-1,k} & {B}_{0,k}
\end{array}
\right),
\end{equation}
where ${B}_{i,k} \in \RR^{p \times p}$, $i=0, \ldots, n-1$. Thus $S_k \in \RR^{N \times N}$ with $N=pn$ has a block-circulant structure. For short we write $S_k:=\Circ(B_{0,k}, \ldots, B_{n-1,k})$.

\begin{remark}\label{rem:norm_Sk}
Due to the structure of $S_k$, it is not difficult to prove that
$$
\|S_k\|_{\infty} \leq \sum_{i=0}^{n-1} \|B_{i,k}\|_{\infty}.
$$
\end{remark}

\smallskip \noindent
If the $k$-th level mesh contains an extraordinary vertex of valence $n$, the refinement rules in its neighborhood involve $pn+1$ points instead of $pn$: $p$ points in each of the $n$ sectors plus the extraordinary vertex. Thus, to construct the local subdivision matrix $S_k$ we first build the matrix
\begin{equation}\label{eq:Sk_block}
\tilde{S}_{k}=\left(\begin{array}{ccccc}
\tilde{\alpha}_{k} & \tilde{\bbeta}_{k}^T & \tilde{\bbeta}_{k}^T & \cdots & \tilde{\bbeta}_{k}^T \\
\tilde{\bgamma}_{k} & \tilde{B}_{0,k} & \tilde{B}_{1,k}  & \cdots & \tilde{B}_{n-1,k}\\
\tilde{\bgamma}_{k} & \tilde{B}_{n-1,k} & \tilde{B}_{0,k} & \cdots & \tilde{B}_{n-2,k}\\
\vdots & \vdots & \vdots  & \ddots & \vdots\\
\tilde{\bgamma}_{k} & \tilde{B}_{1,k}   & \cdots & \tilde{B}_{n-1,k} & \tilde{B}_{0,k}
\end{array}
\right),
\end{equation}
where $\tilde{\alpha}_k \in \RR, \ \tilde{\bbeta}_k, \tilde{\bgamma}_k \in \RR^{p}$ and $\tilde{B}_{i,k} \in \RR^{p \times p}$, $i=0, \ldots, n-1$. Then, following the method shown in \cite[Example 5.14]{PR08}, we transform the matrix $\tilde{S}_k$ in a block-circulant matrix $S_k$ of the form
\begin{equation}\label{eq:matrix_S}
S_k:=\Circ(B_{0,k}, \ldots, B_{n-1,k})\quad \hbox{with} \quad B_{j,k}=\begin{pmatrix} \frac{\tilde{\alpha}_k}{n} & \tilde{\bbeta}_k^T \smallskip\\ \frac{\tilde{\bgamma}_k}{n} & \tilde{B}_{j,k} \end{pmatrix}, \ \ j=0, \ldots, n-1.
\end{equation}
It follows that $S_k \in \RR^{N\times N}$, with $N=n(p+1)$, has a block-circulant structure. Hence, without loss of generality, we can always assume that the local subdivision matrix $S_k$ has a block-circulant structure with blocks of dimension $m \times m$, where $m=p$ if the $k$-th level mesh contains an extraordinary face and $m=p+1$ if it contains an extraordinary vertex.\\

\begin{figure}[!h]
	\centering
	\resizebox{4.5cm}{!}{	\includegraphics{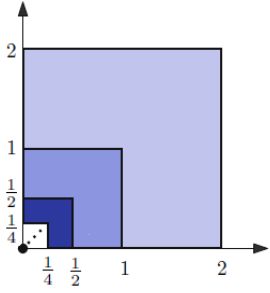}}
	\caption{Domains $\Omega_1, \Omega_2, \Omega_3$ corresponding to three subdivision steps in the case of a quadrilateral mesh containing an extraordinary vertex. (Color figure online.)}
		\label{fig:successione}
\end{figure}

\begin{figure}[h!]
	\centering
	\resizebox{11.0cm}{!}{	\includegraphics{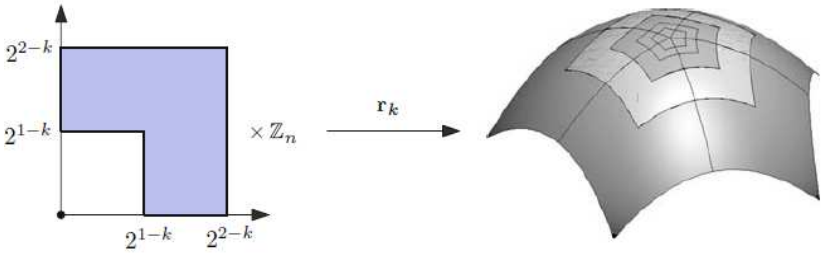}}
	\caption{Ring $\br_k$ in the case of a quadrilateral mesh with an extraordinary vertex (figure taken from \cite{PR08}). (Color figure online.)}
		\label{fig:ring}
\end{figure}

\begin{figure}[!h]
	\centering
	\resizebox{7.5cm}{!}{	\includegraphics{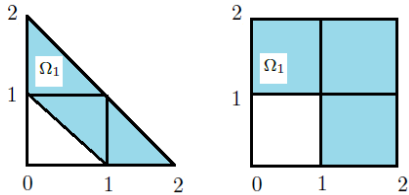}}
	\caption{Domain $\Omega_1$ in the case of a triangular (left) and a quadrilateral (right) mesh containing an extraordinary vertex placed at $\mathbf{0}$. (Color figure online.)}
		\label{fig:dominio}
\end{figure}

\noindent
We continue by introducing some important notation from \cite{PR08,R95,U2000}. We start by assuming that near an isolated extraordinary vertex or face of valence $n$ the  subdivision surface $\br$ is defined on the local domain $\bD_n:=\Omega \times \ZZ_n$ (consisting of $n$ copies of $\Omega$) with
$$\Omega:=\left \{
\begin{array}{ll}
[0,2] \times [0,2] & \hbox{\rm in case of quadrilateral mesh},\\
\{(u,v) \in \RR^2 \, | \, u,v \geq 0 \hbox{ and } 0 \leq u+v \leq 2 \} & \hbox{\rm in case of triangular mesh},
\end{array}
\right.$$
and $\ZZ_n:=\ZZ \slash n\ZZ$.
In the case of triangular and quadrilateral meshes, if we apply one step of refinement to the local domain $\bD_n$, we obtain a new domain with $4n$ cells: $3n$ outer ordinary cells and $n$ inner cells that contain the extraordinary element.
The restriction $\br_1$ of $\br$ to the outer cells is called \emph{ring}.
Denoting by $\tilde{\br}$ the inner part of $\br$, that is $\tilde{\br}:=\br \backslash \br_1$, we can repeat the refinement process only for $\tilde{\br}$ to obtain a second ring $\br_2$ and an even smaller inner part. Hence, iterated refinement generates a sequence of rings $\{\br_k, k \geq 1\}$ which covers all of the surface except for the central point (limit of the sequence of extraordinary vertices or faces), that hereinafter we denote by $\br_c$. Precisely, assuming the central point to be placed at $\mathbf{0}$ and introducing the notation
$$
\tilde{\Omega}:=\left \{
\begin{array}{ll}
[0,1] \times [0,1] & \hbox{\rm in case of quadrilateral mesh},\\
\{(u,v) \in \RR^2 \, | \, u,v \geq 0 \hbox{ and } 0 \leq u+v \leq 1 \} & \hbox{\rm in case of triangular mesh},
\end{array}
\right.
$$
and
$$
\Omega_k:=2^{1-k} (\Omega \backslash \tilde{\Omega}),
\qquad \bD_{n,k}:=\Omega_k \times \ZZ_n, \quad k \geq 1,$$
we see  the ring $\br_k$ as the restriction of the subdivision surface $\br:\bD_n\rightarrow \RR^3$ to the domain $\bD_{n,k}$, i.e.,
$\br_k:=\br |_{\bD_{n,k}}$, (see Figures \ref{fig:successione} and \ref{fig:ring}).
Specifically, in the case of quadrilateral meshes, $\Omega_k$ is explicitly given by
$$\Omega_k= \{(u,v) \in \RR^2 \, | \, u,v \geq 0 \hbox{ and } 2^{1-k} \leq \max\{u,v\} \leq 2^{2-k}\},$$
while in the case of  triangular meshes
$$\Omega_k= \{(u,v) \in \RR^2 \, | \, u,v \geq 0 \hbox{ and } 2^{1-k} \leq u+v \leq 2^{2-k}\},$$
(see Figure \ref{fig:dominio}). As a consequence, both in the case of triangular and quadrilateral meshes, $\Omega_k$ is constituted by the union of 3 cells, say $\omega_k^{[1]}$, $\omega_k^{[2]}$ and $\omega_k^{[3]}$, implying  that the domain $\bD_{n,k}$ is indeed made of $3n$ cells. It follows that the entire surface ring $\br_k$ is the union of $3n$ patches, each one denoted by $\br_k^{[j]}$ and corresponding to the restriction of the subdivision surface $\br$ to the single cell $\omega_k^{[j]}$, $j \in \mathbb{J}_{3n}$ where
\begin{equation}\label{eq:J3n}
\mathbb{J}_{3n}:=\{j\in \ZZ \ : \ j=3(l-1)+i, \ \ l=1,\ldots,n, \ \ i=1,2,3\}.
\end{equation}

Exploiting the given definition of $\br_{k}$, we can now provide
the following notion of convergence of a non-stationary subdivision scheme ${\mathscr S}$  in irregular regions.

\begin{definition}\label{def:conv}
Let ${\mathscr S}$ be a (non-stationary) subdivision scheme with the property of convergence in regular regions and whose action in an irregular region is described by a matrix sequence  $\{S_k \in \RR^{N \times N}, \, k \geq 1\}$. Moreover, let $\bd_1 \in \RR^{N \times 3}$ be the vector with the vertices of $\cE^{(1)}$.
${\mathscr S}$ is said to be convergent in the neighborhood of an extraordinary vertex/face of valence $n$ if, for all initial data $\bd_1$, there exists a limit point $\br_c \in\mathbb{R}^3$ such that
\begin{equation}\label{eq:conv_def}
\lim_{k\rightarrow +\infty} \sup_{(u,v) \in \bD_{n,k}} \| \br_{k}(u,v)-\br_c\|_{\infty}=0.
\end{equation}
\end{definition}

\noindent
We conclude by observing that, if the subdivision scheme ${\mathscr S}$ converges then $\br= \bigcup_{k \geq 1} \br_k \cup \{ \br_c \}$ is a surface without gap, i.e., $\br$ is a surface which is continuous at all points including $\br_c$.
The surface $\br$ is called the \emph{limit surface} of the  subdivision scheme ${\mathscr S}$.

\medskip
In the following, for a subdivision scheme ${\mathscr S}$ with the property of convergence in regular regions, we denote by $\bd_k^{[j]} \in \RR^{P\times 3}$, $P<N$,  the vector with the control points of each patch $\br_k^{[j]}$, and with $\bPhi_k^{[j]}\in \RR^{P}$ the function vector containing all the basic limit functions $\phi_k$ of ${\mathscr S}$ whose supports intersect $\omega_k^{[j]}$. We assume that the functions in $\bPhi_k^{[j]}$ are ordered as the points in the vector $\bd_k^{[j]}$, and thus we call them the associated basic limit functions.
Therefore, for each $j \in \mathbb{J}_{3n}$, we have
\begin{equation}\label{eq:r_k}
\begin{array}{c}
\displaystyle \br_k^{[j]}: \omega_k^{[j]} \rightarrow \RR^3 \smallskip\\
\displaystyle (u,v)\longmapsto \br_k^{[j]} (u,v)=(\bd_k^{[j]})^T \, \bPhi_k^{[j]}(u,v).
\end{array}
\end{equation}

\smallskip
Now, assume also that $\bar{\mathscr S}$ is a stationary subdivision scheme with the property of convergence in regular regions.
Denoting by $\bar{\bPhi}^{[j]} \in \RR^{P}$ the vector containing all the basic limit functions $\bar{\phi}$ of $\bar{\mathscr S}$ whose supports intersect $\omega_k^{[j]}$, if the assumptions of Theorem \ref{prop:asympt_equiv} are satisfied, we have that
$$
\lim_{k\rightarrow +\infty} \sup_{(u,v) \in \omega_{k}^{[j]}} \|\bPhi_k^{[j]}(u,v)-\bar{\bPhi}^{[j]}(u,v)\|_{\infty}=0, \quad \forall j \in \mathbb{J}_{3n}.
$$

\smallskip
\begin{remark}\label{rem:POU}
Let $\bx_0:=(1, 1, \ldots, 1)^T \in \RR^P$. We observe that, for a convergent stationary subdivision scheme $\bar{\mathscr S}$,  we have  $\overline{\bPhi}^{[j]}(u,v)^T \bx_0=1$ for all $(u,v)\in\RR^2$, $j \in \mathbb{J}_{3n}$. Instead, for a  non-stationary subdivision scheme ${\mathscr S}$ with the property of convergence in regular regions, ${\bPhi}_k^{[j]}(u,v)^T \bx_0=1$ for all $k\geq 1$ and  for all $(u,v)\in\RR^2$, $j \in \mathbb{J}_{3n}$ if and only if ${\mathscr S}$ has the property of stepwise reproduction of constants (see, e.g., \cite{CCR14} for more details). In general, ${\bPhi}_k^{[j]}(u,v)^T \bx_0=\alpha_k$ with $\alpha_k \in \RR$, for all $j \in \mathbb{J}_{3n}$.
\end{remark}
\noindent
Now, let $\bd_1 \in \RR^{N \times 3}$ be the collection of the vectors of control points $\bd_1^{[j]}$ of all patches $\br_1^{[j]}$, $j \in \mathbb{J}_{3n}$. Denoted by $\{S_k \in \RR^{N \times N}, k \geq 1\}$ the matrix sequence that defines
a non-stationary subdivision scheme ${\mathscr S}$ in an irregular region, we can obtain
the entire set of the $(k+1)$-th level control points representing the whole ring $\br_{k+1}$ by the matrix multiplication
\begin{equation}\label{eq:defS(k)}
\bd_{k+1}=S_k \bd_k=S_k S_{k-1} \bd_{k-1}=...=S^{(k)} \bd_1 \; \hbox{with} \; S^{(k)}:= \left \{
\begin{array}{ll}
S_{k} \, S_{k-1} \, \cdots \, S_1, & k \geq 1,\\
I, & k=0.
\end{array}
\right .
\end{equation}
Moreover, denoting by ${\bPhi}_{k+1}$ the function vector with blocks ${\bPhi}_{k+1}^{[j]}$, $j \in \mathbb{J}_{3n}$, we can rewrite each patch $\br_{k+1}^{[j]} (u,v)=(\bd_{k+1}^{[j]})^T \, \bPhi_{k+1}^{[j]}(u,v)$ of the surface ring $\br_{k+1}$ as
$$\br_{k+1}^{[j]} (u,v)=\bd_{k+1}^T \, \bPhi_{k+1}(u,v), \qquad (u,v)\in \omega_{k+1}^{[j]},$$
(\ie, independently of $j$) since the function vector $\bPhi_{k+1} \in \RR^{N}$ indeed contains only $P$ functions that are non-zero on $\omega_{k+1}^{[j]}$.

\bigskip
The goal of the next section is to prove new basic results that allow us to derive a general criterion for verifying if the limit surface $\br$ generated by a rotationally symmetric, non-stationary subdivision scheme is normal continuous.

\begin{definition}\label{def:G1}
The surface $\br$, limit of a convergent non-stationary subdivision scheme ${\mathscr S}$ which is $C^1$-convergent in regular regions, is \emph{normal continuous} at $\br_c$ (limit point of a sequence of extraordinary vertices/faces of valence $n$) if there exists a unique vector $\mathbf{n}_\infty$ such that
$$
\lim_{k\rightarrow +\infty} \sup_{(u,v) \in \bD_{n,k}} \| \bn_{k}(u,v)-\mathbf{n}_\infty\|_{\infty}=0,
$$
for almost all sequences of normal vectors $\{\bn_{k}(u,v):=\frac{ \partial_u \br_{k}(u,v) \wedge \partial_v \br_{k}(u,v)}{\|\partial_u \br_{k}(u,v) \wedge \partial_v \br_{k}(u,v)\|_2}, \, k \geq 1\}$,
where $\br_{k}(u,v)$ satisfies \eqref{eq:conv_def}.
\end{definition}

\section{New results linked to the $C^1$-convergence analysis of non-stationary subdivision schemes in regular regions}\label{sec:new3}

The preliminary results required in Subsection \ref{sec:4} to give sufficient conditions
for verifying normal continuity of the subdivision surface deal with new results connected with the
$C^1$-convergence analysis of non-stationary subdivision schemes in regular regions.
For them we recall the well-known notions of asymptotical equivalence of order $1$ and of divided-difference scheme, plus related results proven in~\cite{DL95}.

\begin{definition}
Let ${\mathscr S}$ and $\bar{\mathscr S}$ be subdivision schemes defined in regular regions by the subdivision masks $\{\bc^{(k)}\in \ell(\ZZ^2), \, k \geq 1\}$ and $\bc \in \ell(\ZZ^2)$, respectively.
If
$$
\sum_{k=1}^{+\infty} 2^k \|{\cal S}_{\bc^{(k)}}-{\cal S}_{\bc}\|_{\infty} < +\infty,
$$
then ${\mathscr S}$ and $\bar{\mathscr S}$ are said to be \emph{asymptotically equivalent schemes of order $1$}.
\end{definition}

\begin{remark}
Asymptotical equivalence of order $1$ implies asymptotical equivalence in the sense of Definition \ref{def:asympt_eq}.
\end{remark}

\begin{theorem}\cite[Theorem 8]{DL95}\label{prop:asympt_equiv_order1}
Let ${\mathscr S}$ and $\bar{\mathscr S}$ be subdivision schemes defined in regular regions by the subdivision masks $\{\bc^{(k)} \in \ell(\ZZ^2), \, k \geq 1\}$ and $\bc \in \ell(\ZZ^2)$, respectively.
If ${\mathscr S}$ and $\bar{\mathscr S}$ are asymptotically equivalent of order $1$,
then
$C^1$-convergence of $\bar{\mathscr S}$ implies $C^1$-convergence of ${\mathscr S}$.
\end{theorem}

\begin{definition}
For the two perpendicular directions $\be_1=(1,0)^T$, $\be_2=(0,1)^T$, we define as
$$(\Delta_{\be_j}^{(\ell)} \bbf^{(\ell)})_\balpha:=\frac{\rf^{(\ell)}_\balpha-\rf^{(\ell)}_{\balpha-\be_j}}{2^{-\ell}}, \quad \balpha \in \ZZ^2, \ j \in \{1,2\}, \ \ell \geq 1,$$
the ${\be}_j$-directional divided difference operator.
\end{definition}

The following lemma recalls a well-known property fulfilled by the symbols of the so-called divided difference schemes. Its proof is omitted since already given in \cite[Section 4.2.2]{DL02}.

\begin{lemma}\label{lem:divdiff_scheme}
Let $j \in \{1,2\}$. If $c^{(\ell)}(\bz)=\frac12 (1+z_j) b_{\be_j}^{(\ell)}(\bz)$, then
$$
\Delta_{\be_j}^{(\ell+1)} f^{(\ell+1)}(\bz)=b_{\be_j}^{(\ell)}(\bz) \ \Delta_{\be_j}^{(\ell)} f^{(\ell)}(\bz^2),
$$
and $\{{\cal S}_{\bb_{\be_j}^{(\ell)}}, \, \ell\geq 1 \}$ is called the $\be_j$-directional divided difference scheme of $\{{\cal S}_{\bc^{(\ell)}}, \, \ell\geq 1 \}$.
\end{lemma}

\smallskip \noindent
From Lemma \ref{lem:divdiff_scheme} we have that
\begin{equation}\label{def:Sb_scheme}
\Delta_{\be_j}^{(\ell+1)} \bbf^{(\ell+1)}={\cal S}_{\bb_{\be_j}^{(\ell)}} \, \Delta_{\be_j}^{(\ell)} \bbf^{(\ell)}
\quad \Leftrightarrow \quad
\Delta_{\be_j}^{(\ell+1)} {\cal S}_{\bc^{(\ell)}} \, \bbf^{(\ell)}={\cal S}_{\bb_{\be_j}^{(\ell)}} \, \Delta_{\be_j}^{(\ell)} \bbf^{(\ell)}.
\end{equation}

\begin{lemma}\label{lem:AE_inherit}
Let ${\mathscr S}$ and $\bar{\mathscr S}$ be subdivision schemes specified in regular regions by the subdivision symbols $\{c^{(k)}(\bz), \, k \geq 1\}$ and $c(\bz)$, respectively. Assume that:
\begin{itemize}
\item[i)] ${\mathscr S}$ and $\bar{\mathscr S}$ are asymptotically equivalent of order $1$;
\item[ii)] the factor $(1+z_1)(1+z_2)$ is contained in the symbols $c(\bz)$ and $c^{(k)}(\bz)$, for all $k \geq 1$.
\end{itemize}
Then, the divided difference schemes with  symbols
$b_{\be_j}(\bz):=\frac{2 c(\bz)}{1+z_j}$, $j\in \{1,2\}$ and $b_{\be_j}^{(k)}(\bz):=\frac{2 c^{(k)}(\bz)}{1+z_j}$, $j\in \{1,2\}$,
are asymptotically equivalent of order $1$.
\end{lemma}

\proof
We only consider the case corresponding to $j=1$, since the case $j=2$ can be treated analogously.
To simplify the notation we denote $b_{\be_1}(\bz)$  and $b^{(k)}_{\be_1}(\bz)$ by $b(\bz)$  and $b^{(k)}(\bz)$, respectively. We start by considering the relation
$$2c(\bz)=(1+z_1)b(\bz)$$
with
$$c(\bz):=\sum_{\balpha\in [0,N_1]\times[0,N_2]} {\rm c}_{\balpha}\bz^{\balpha} \qquad \hbox{and} \qquad b(\bz):=\sum_{\balpha\in [0,N_1-1]\times[0,N_2]} {\rm b}_{\balpha}\bz^{\balpha}.$$
Comparing the same power of $\bz$ we easily see that,
$$
{\rm c}_{0,\alpha_2}=\frac{1}{2} {\rm b}_{0,\alpha_2}, \;
{\rm c}_{N_1,\alpha_2}=\frac{1}{2} {\rm b}_{N_1-1,\alpha_2},
\quad {\rm c}_\balpha=\frac12 \left( {\rm b}_{\balpha}+ {\rm b}_{\balpha-\be_1} \right),\; \balpha\in [1,N_1-1]\times[0,N_2],
$$
which means
$$
{\small
{\rm b}_{0,\alpha_2}=2 {\rm c}_{0,\alpha_2},\; {\rm b}_{N_1-1,\alpha_2}=2 {\rm c}_{N_1,\alpha_2}, \;  {\rm b}_{\balpha}=2\sum_{\beta_1=0}^{\alpha_1}(-1)^{\alpha_1-\beta_1} {\rm c}_{\beta_1,\alpha_2},\;  \balpha\in [1,N_1-2]\times[0,N_2]
}.$$
Analogously, working with the relation $2c^{(k)}(\bz)= (1+z_1)b^{(k)}(\bz)$ we get
$$
{\small
{\rm b}^{(k)}_{0,\alpha_2}=2 {\rm c}^{(k)}_{0,\alpha_2},\, {\rm b}^{(k)}_{N_1-1,\alpha_2}=2 {\rm c}^{(k)}_{N_1,\alpha_2},\,  {\rm b}^{(k)}_{\balpha}=2\sum_{\beta_1=0}^{\alpha_1}(-1)^{\alpha_1-\beta_1} {\rm c}^{(k)}_{\beta_1,\alpha_2},\; \balpha\in [1,N_1-2]\times[0,N_2]
}.$$
Therefore,
$$\| \bb^{(k)}-\bb\|_{\infty} \leq  2N_1  \| \bc^{(k)}-\bc\|_{\infty} \qquad \hbox{and} \qquad  \sum_{k=1}^{+\infty} 2^{k}\| \bb^{(k)}-\bb\|_{\infty} \leq  2N_1  \sum_{k=1}^{+\infty} 2^{k}\| \bc^{(k)}-\bc\|_{\infty} <\infty.
$$
Thus, in light of Remark \ref{rem:CDMM}, the result is proven.
\endproof

\smallskip \noindent
The previous Lemma is useful for the next result.

\begin{proposition}\label{prop:conv_Sbscheme}
Let ${\mathscr S}$ and $\bar{\mathscr S}$ be subdivision schemes such that in regular regions:
\begin{itemize}
\item[i)] ${\mathscr S}$ and $\bar{\mathscr S}$ are asymptotically equivalent of order $1$.
\item[ii)] $\bar{\mathscr S}$ is $C^1$-convergent with symbol $c(\bz)$ that contains the factor $(1+z_1)(1+z_2)$;
    \item[iii)]${\mathscr S}$ is defined by the subdivision symbols $\{ c^{(\ell)}(\bz), \, \ell \geq 1\}$
all containing the factor $(1+z_1)(1+z_2)$.
\end{itemize}
Then, the associated divided difference schemes with symbols $b_{\be_j}(\bz):=\frac{2 c(\bz)}{1+z_j}$ and $b_{\be_j}^{(\ell)}(\bz):=\frac{2 c^{(\ell)}(\bz)}{1+z_j}$, $j\in \{1,2\}$, satisfy the following properties:
\begin{itemize}
\item[a)] the sequence of basic limit functions of $\{{\cal S}_{\bb_{\be_j}^{(\ell)}}, \, \ell \geq 1 \}$
converges uniformly to the basic limit function of $\{{\cal S}_{\bb_{\be_j}}\}$;
\item[b)]
$
\displaystyle{\lim_{\ell \rightarrow +\infty} \, {\cal S}_{\bb_{\be_j}^{(k+\ell)}} \, {\cal S}_{\bb_{\be_j}^{(k+\ell-1)}} \, \ldots {\cal S}_{\bb_{\be_j}^{(k)}} \,  \Delta_{\be_j}^{(k)} \bdelta=\partial_{\be_j} \phi_k}
\quad \hbox{and} \quad
\lim_{\ell \rightarrow +\infty} \, ({\cal S}_{\bb_{\be_j}})^{\ell} \Delta_{\be_j}^{(1)} \bdelta=\partial_{\be_j} \overline{\phi},
$\\
for $\bdelta=\{\delta_{0,\balpha}, \, \balpha \in \ZZ^2\}$ and with $\phi_k$ defined as in \eqref{def:BLFnonstat} and $\overline{\phi}$ as in \eqref{def:BLFstat}.
\end{itemize}

\end{proposition}

\proof
The result in {\it a}) is a direct consequence of Lemma \ref{lem:AE_inherit} and Theorem \ref{prop:asympt_equiv} (see also \cite[Lemma 15]{DL95}).\\
To show {\it b}) we proceed as follows. In view of the factorization properties of $c^{(\ell)}(\bz)$, we can apply Lemma \ref{lem:divdiff_scheme} to conclude the existence of the $\be_j$-directional divided difference scheme of order 1 of $\{{\cal S}_{\bc^{(\ell)}}, \, \ell\geq 1 \}$.
Then, to show convergence of the $\be_j$-directional divided difference scheme of order 1, we just recall the result in {\it a}).
Next, we exploit \eqref{def:Sb_scheme} and write
$$
{\cal S}_{\bb_{\be_j}^{(k)}} \, \Delta_{\be_j}^{(k)} \bdelta=\Delta_{\be_j}^{(k+1)} {\cal S}_{\bc^{(k)}} \, \bdelta,
$$
so that
$$
{\cal S}_{\bb_{\be_j}^{(k+\ell)}} \, {\cal S}_{\bb_{\be_j}^{(k+\ell-1)}} \, \ldots {\cal S}_{\bb_{\be_j}^{(k)}} \,  \Delta_{\be_j}^{(k)} \bdelta=
\Delta_{\be_j}^{(k+\ell+1)} {\cal S}_{\bc^{(k+\ell)}} \, {\cal S}_{\bc^{(k+\ell-1)}} \, \ldots \,  {\cal S}_{\bc^{(k)}} \, \bdelta.
$$
Moreover, introducing the notation $\bdelta^{(k+\ell+1)}:={\cal S}_{\bc^{(k+\ell)}} \, {\cal S}_{\bc^{(k+\ell-1)}} \, \ldots \,  {\cal S}_{\bc^{(k)}} \, \bdelta$, we have that
$$
\Delta_{\be_j}^{(k+\ell+1)} \bdelta^{(k+\ell+1)}=\frac{\bdelta^{(k+\ell+1)}- (\bdelta^{(k+\ell+1)})_{\cdot-\be_j}}{2^{-k-\ell-1}}, \quad j \in \{1,2\}.
$$
Thus
$$
\begin{array}{l}
\lim_{\ell \rightarrow +\infty} \, {\cal S}_{\bb_{\be_j}^{(k+\ell)}} \, {\cal S}_{\bb_{\be_j}^{(k+\ell-1)}} \, \ldots {\cal S}_{\bb_{\be_j}^{(k)}} \,  \Delta_{\be_j}^{(k)} \bdelta=\\
\lim_{\ell \rightarrow +\infty} \,  \Delta_{\be_j}^{(k+\ell+1)} {\cal S}_{\bc^{(k+\ell)}} \, {\cal S}_{\bc^{(k+\ell-1)}} \, \ldots \,  {\cal S}_{\bc^{(k)}} \, \bdelta =\\
\lim_{\ell \rightarrow +\infty} \, \frac{{\cal S}_{\bc^{(k+\ell)}} \, {\cal S}_{\bc^{(k+\ell-1)}} \, \ldots \,  {\cal S}_{\bc^{(k)}} \, \bdelta -  \big( {\cal S}_{\bc^{(k+\ell)}} \, {\cal S}_{\bc^{(k+\ell-1)}} \, \ldots \,  {\cal S}_{\bc^{(k)}} \, \bdelta \big)_{\cdot-\be_j}}{2^{-k-\ell-1}} =\\
\partial_{\be_j} \phi_k,
\end{array}
$$
in view of the fact that $\lim_{\ell \rightarrow +\infty}  {\cal S}_{\bc^{(k+\ell)}} \, {\cal S}_{\bc^{(k+\ell-1)}} \, \ldots {\cal S}_{\bc^{(k)}} \,  \bdelta=\phi_k$ and $\phi_k$ is $C^1$.\\
The result for the stationary scheme follows by taking ${\cal S}_{\bc^{(\ell)}}={\cal S}_{\bc}$ for all $\ell \geq 1$ and using Theorem \ref{prop:asympt_equiv}.
\endproof

\smallskip \noindent
As a consequence of the  previous proposition we have
\begin{corollary}\label{cor:unif_conv_der}
 Under the assumptions of Proposition \ref{prop:conv_Sbscheme}
$$
\lim_{k\rightarrow + \infty} \sup_{(u,v) \in \RR^2} | \partial_{\be_j} \phi_k(u,v) - \partial_{\be_j} \overline{\phi}(u,v) |=0, \qquad j \in \{1,2\}.
$$
\end{corollary}

\section{Analysis of rotationally symmetric, non-stationary subdivision schemes in irregular regions}\label{sec:3}

Before focusing on the sufficient conditions that guarantee the convergence of a rotationally symmetric, non-stationary subdivision scheme in the neighborhood of an extraordinary element (Theorem \ref{theo:convergenceNEW}), we present a few linear algebra results to be used for the subdivision analysis.

\subsection{Auxiliary linear algebra results}\label{subsec:algebra}
Let $M \in \RR^{N \times N}$.
In the following, two simple results based on the Jordan decomposition of $M$ are proven.
For the first one we assume $\bd\in \RR^{N\times 3}$ and consider the sequence $\{M^k \bd, \, k\geq 0\}$. Then, under suitable assumptions on the matrix $M$, we show its convergence.
 In the second one (which is a well known result so that we omit its proof) we study the properties of $M^k$, $k \geq 0$, again with the help of its Jordan decomposition.

\begin{proposition}\label{prop:conv_stat}
Assume that $M$ has the unique dominant eigenvalue $1$ and that the associated eigenvector is $\bx_0=(1,1,...,1)^T$. Let $XJX^{-1}$ be the Jordan decomposition of $M\in \RR^{N\times N}$ and let $\bx_0$ be the first column of $X$.
Then, for all $\bd\in \RR^{N\times 3}$,
\begin{equation}\label{eq:lim}
\lim_{k\rightarrow +\infty} {M}^k\bd= \bx_0 \bq^T,
\end{equation}
with $\bq^T= \tilde{\bx_0}^T \bd\in \RR^{1\times 3}$, $\tilde{\bx_0}^T=\be_1^TX^{-1}\in \RR^{1\times N}$ and  $\be_1^T=(1,0,...,0)\in \RR^{1\times N}$. Moreover,  $\tilde{\bx_0}^T M= \tilde{\bx_0}^T$.
\end{proposition}

\proof
Using the Jordan decomposition of $M$ we can write $M^k=X J^k X^{-1}.$ Hence, recalling that $1$ is the unique dominant eigenvalue of $J$ and the associated eigenvector is $\bx_0=(1,1,...,1)^T$, we have
$$
\lim_{k\rightarrow + \infty} J^k=\begin{pmatrix} 1 & \cdots&  0 & 0\\
                0 & 0 & \cdots & 0\\
								\vdots & \ddots & \ddots & \vdots\\
								0 & \cdots & \cdots & 0 \end{pmatrix}=\be_1 \be_1^T \qquad
\hbox{and}
\qquad
\be_1^T J =\be_1^T.
$$
Therefore,
$$
\lim_{k\rightarrow +\infty} M^k \bd=
\lim_{k\rightarrow +\infty }X J^k X^{-1}\bd=
X \left(\lim_{k\rightarrow + \infty} J^k \right)X^{-1}\bd=
X \be_1 \be_1^T X^{-1}\bd,
$$
and, in view of the fact that $X \be_1=\bx_0$, \eqref{eq:lim} is proven.
Moreover,
$$\tilde{\bx_0}^T M=\be_1^T X^{-1}M=\be_1^T X^{-1}X J X^{-1}=\be_1^T J X^{-1}=\be_1^T X^{-1}=\tilde{\bx_0}^T,$$
so concluding the proof.
\endproof

\noindent
In Propositions \ref{prop:M_stable} and \ref{lemmanew}, $\| \cdot \|$ refers to any vector norm and its induced matrix norm.

\begin{proposition}\label{prop:M_stable}
If the dominant eigenvalue of $M$ is $1$ and its algebraic multiplicity  is $1$, then there exists a finite positive constant $\mathpzc{C}$ (independent of $k$) such that
$$\|M^k\| \leq \mathpzc{C}, \quad \forall k \geq 0.$$
\end{proposition}

\smallskip
\begin{remark}\label{rem:convergenceimplies}
It is important to remark that, according to \cite[Theorem 4.20]{PR08},
we can assume without loss of generality that the subdivision matrix $S$ defining a rotationally symmetric, stationary subdivision scheme $\bar{\mathscr S}$ does not have ineffective eigenvectors.
Thus, hereafter $S$ always satisfies the assumptions of Propositions \ref{prop:conv_stat} and \ref{prop:M_stable} since (see, e.g., \cite{PR08,R95,Z00})
\begin{itemize}
\item[1)] the unique dominant eigenvalue of $S$ is $\lambda_0=1$,
\item[2)] the algebraic multiplicity of $\lambda_0$ is $1$,
\item[3)] the eigenvector associated with $\lambda_0$ is $\bx_0=(1,1,...,1)^T$.
\end{itemize}
\end{remark}

\smallskip \noindent
In the next Proposition we replace the $k$-th power of the matrix $M$ with the product of $k$ different matrices $M_{k} \, M_{k-1} \, \cdots \, M_1$ and we successively consider hybrid combinations of the two.

\begin{proposition} \label{lemmanew}
Let $M^{(0)}:=I \in \RR^{N\times N}$ and let $M^{(k)}:=M_{k} \, M_{k-1} \, \cdots \, M_1$ with $M_j \in \RR^{N \times N}$, for all $j=1,...,k$ and for all $k \geq 1$.
Let $M\in \RR^{N\times N}$ be a nonsingular matrix having $1$ as dominant eigenvalue with algebraic multiplicity $1$.
If, for all $k \geq 1$, $\|M_k-M\| \leq \frac{\mathpzc{C}}{\sigma^k}$ with $\sigma>1$ and some finite positive constant $\mathpzc{C}$ (independent of $k$), then
$$\|M^{(k)}\| \leq \mathpzc{\widehat C}, \quad  \forall \ k \geq 1,$$
with $\mathpzc{\widehat C}$ a finite positive constant (independent of $k$).
\end{proposition}
\proof
The proof takes inspiration from \cite[Theorem 5]{DL95}. The claim is proven by introducing, for $\by \in \RR^N$, a new vector norm
$$
\|\by\|_{M}:=\sup_{k \geq 0} \|M^{k} \by\|,
$$
associated to the given nonsingular matrix $M\in \RR^{N\times N}$.
In view of Proposition \ref{prop:M_stable}, our assumption on $M$ implies the existence of a finite positive constant $\widetilde{\mathpzc{C}}$ such that $\|M^k\| \leq \widetilde{\mathpzc{C}}$ for all $k \geq 0$.
Moreover, $\|\by\| \leq \|\by\|_{M}$ since $\|M^{k} \by\|=\|\by\|$ when $k=0$.
There follows that
$$
\|\by\| \leq \|\by\|_{M}\leq \widetilde{\mathpzc{C}} \|\by\|,
$$
meaning that any standard vector norm and the $\| \cdot \|_{M}$ norm are uniformly equivalent.
We now consider the induced norm for the matrix $M$ itself, and denote it as $\|\cdot\|_{M}$.
Then
$$
\begin{array}{lll}
\|M\|_{M}&:=& \displaystyle \sup_{\|\by\|_{M}=1} \|M \by\|_{M}
=\sup_{\|\by\|_{M}=1} \, \sup_{k \geq 0} \|M^{k}(M\by)\|\\
&=& \displaystyle \sup_{\|\by\|_{M}=1} \, \sup_{k \geq 1} \|M^{k}\by\| \leq \sup_{\|\by\|_{M}=1} \, \sup_{k \geq 0} \|M^{k}\by\|=\sup_{\|\by\|_{M}=1} \|\by\|_{M}=1.
\end{array}
$$
We continue by exploiting the uniform equivalence of norms to bound $\|M^{(k)}\|$, $k \geq 1$. Obviously,
$$
\|M_k M_{k-1} ... M_1 \|_{M}
\leq \|M_k\|_{M} \, \|M_{k-1}\|_{M} \, ... \, \|M_{1}\|_{M} \, .
$$
Furthermore,
$$
\|M_j-M\|_{M}\leq \widetilde{\mathpzc{C}} \|M_j-M\| \leq \widetilde{\mathpzc{C}}  \frac{\mathpzc{C}}{\sigma^j},
$$
where $\mathpzc{C}$ is the finite positive constant appearing in the assumption and where $\widetilde{\mathpzc{C}}$, different from above,  comes from the norm equivalence.
Finally, for any arbitrary $k \geq 1$, we arrive at:
$$
\begin{array}{lll}
\|M^{(k)}\| &\leq& \prod_{j=1}^k \|M_j\|_{M} \leq \displaystyle \, \prod_{j=1}^k \left( \|M\|_{M} + \|M_j-M\|_{M} \right )\smallskip \\
&\leq& \displaystyle \, \prod_{j=1}^k \left( 1 +   \frac{\widetilde{\mathpzc{C}} \mathpzc{C}}{\sigma^j} \right )
= \displaystyle \e^{\log_e \left( \prod_{j=1}^k \big( 1+\frac{\widetilde{\mathpzc{C}} \mathpzc{C}}{\sigma^j} \big)\right)}
= \displaystyle \e^{\sum_{j=1}^k  \log_e  \big( 1+\frac{\widetilde{\mathpzc{C}} \mathpzc{C}}{\sigma^j} \big) }\\
&\leq& \displaystyle \e^{\sum_{j=1}^k  \frac{\widetilde{\mathpzc{C}} \mathpzc{C}}{\sigma^j} }
\leq \displaystyle \e^{ \widetilde{\mathpzc{C}} \mathpzc{C} \sum_{j=1}^{+\infty}  \frac{1}{\sigma^j} },
\end{array}
$$
where the last but one inequality follows from the fact that $\log_e(1+x) \leq  x$ for all $x \geq 0$.
Since  $\sum_{j=1}^{+\infty} \frac{1}{\sigma^j}<\infty$ the claim follows.
\endproof

\smallskip \noindent
We conclude this section with another useful intermediate result relating ${M}^{(k)}$ with $M^k$.

\begin{proposition}\label{prop:Mk_ricorrenza}
Let $M \in \RR^{N \times N}$, $M^{(0)}:=I\in \RR^{N \times N}$,  $M^{(k)}:=M_k\cdots M_1 \in \RR^{N \times N}$, for all $k \geq 1$. Then,
$$
 {M}^{(k)}= {M}^k+\sum_{j=1}^{k} {M}^{k-j}(M_j-M)  {M}^{(j-1)},\quad \hbox{for all}\quad k \geq 1.
$$
\end{proposition}
\proof
Assuming $\sum_{j=1}^{k-1} {M}^{k-j}(M_j-M){M}^{(j-1)}$ to be $\mathbf{0}$ when $k=1$, we can write
$$\begin{array}{l}
\displaystyle {M}^k+\sum_{j=1}^{k} {M}^{k-j}(M_j-M)  {M}^{(j-1)}=\\
\displaystyle {M}^k+\sum_{j=1}^{k-1} {M}^{k-j}(M_j-M){M}^{(j-1)}+(M_k-M){M}^{(k-1)}=\\
\displaystyle {M}^k+M^{(k)}+\sum_{j=1}^{k-1} {M}^{k-j}M^{(j)}-MM^{(k-1)}-\sum_{j=1}^{k-1}M^{k-j+1}M^{(j-1)}=\\
\displaystyle M^k+M^{(k)}+\sum_{j=1}^{k-2}M^{k-j}M^{(j)}-\sum_{j=0}^{k-2}M^{k-j}M^{(j)}=\\
\displaystyle M^{k}+M^{(k)}-M^{k}=
\displaystyle M^{(k)},
\end{array}$$
so concluding the proof.
\endproof

\subsection{Convergence analysis in irregular regions}\label{subsec:convergence}
In this section we make use of the previous linear algebra results to provide sufficient conditions for establishing the convergence of a rotationally symmetric, non-stationary subdivision scheme ${\mathscr S}$ defined in an irregular region by a matrix sequence $\{S_k \in \RR^{N \times N}, k \geq 1\}$.
With the notation previously introduced, let $\bd_{k+1} \in \RR^{N \times 3}$ be the collection of the vectors of control points $\bd_{k+1}^{[j]}$ of all patches $\br_{k+1}^{[j]}$, $j \in \mathbb{J}_{3n}$ with $\mathbb{J}_{3n}$ given in \eqref{eq:J3n}.
According to \eqref{eq:defS(k)}, the entire set of the $(k+1)$-th level control points $\bd_{k+1}$ representing the whole ring $\br_{k+1}$ is given by the matrix multiplication
$$
\bd_{k+1}=S^{(k)} \bd_1 \quad \hbox{with} \quad S^{(k)}:= \left \{
\begin{array}{ll}
S_{k} \, S_{k-1} \, \cdots \, S_1, & k \geq 1,\\
I, & k=0.
\end{array}
\right .
$$
Recalling Definition \ref{def:conv}, our goal is to study the convergence of the sequence of regular rings $\{\br_{k+1}, k \geq 0\}$ whose patches $\br_{k+1}^{[j]}$ are described by the equation
$$\br_{k+1}^{[j]}=(\bd_{k+1}^{[j]})^T \, \bPhi_{k+1}^{[j]}=\bd_{k+1}^T \, \bPhi_{k+1}, \quad j \in \mathbb{J}_{3n}.$$

\noindent
The key idea to prove convergence of ${\mathscr S}$ is to write the product matrix $S^{(k)}$ in terms of the stationary matrix $S^k$. Indeed, from Proposition \ref{prop:Mk_ricorrenza} we write
\begin{equation}\label{eq:global_cpset}
\bd_{k+1}={S}^k \bd_1 + \by_k
\quad \hbox{with} \quad
\by_k:=\sum_{j=1}^{k} {S}^{k-j}(S_j-S) {S}^{(j-1)}\bd_1,
\end{equation}
and then show our first main result.

\begin{theorem}\label{theo:convergenceNEW}
Let ${\mathscr S}$ be a rotationally symmetric, non-singular, non-stationary subdivision scheme whose action in an irregular region is described by a matrix sequence  $\{S_k, k \geq 1\}$.
Moreover, let $\bar{\mathscr S}$ be a symmetric, stationary subdivision scheme that in the same irregular region is associated with $S$.
Assume that:
\begin{itemize}
\item[(i)] $\bar{\mathscr S}$  is convergent both in regular and  irregular regions,

\item[(ii)] ${\mathscr S}$ is asymptotically equivalent to $\bar{\mathscr S}$ in regular regions,

\item[(iii)] in the irregular region the matrices $S_k$ and $S$ satisfy, for all $k \geq 1$,  $\|S_k-S\|_{\infty} \leq \frac{\mathpzc{C}}{\sigma^k}$ with $\mathpzc{C}$ some finite positive constant and $\sigma>1$.
\end{itemize}
Then, for all initial data $\bd_1 \in \RR^{N \times 3}$, the non-stationary subdivision scheme ${\mathscr S}$ is convergent also in the irregular region. In particular,
$$
\lim_{k \rightarrow + \infty} \sup_{(u,v) \in \Omega_{k+1}} \|\br_{k+1}(u,v)-(\bq_0+\bbeta_0)\|_{\infty}=0,
$$
where
\begin{itemize}
\item $\bq_0=\bd_1^T \tilde{\bx}_0 \in \RR^3$ with $\tilde{\bx}_0$ such that $S^T \tilde{\bx}_0=\tilde{\bx}_0$,
\\
\item $\bbeta_0= (\lim_{k\rightarrow +\infty} \by_k)^T \, \frac{\bx_0}{N}   \in \RR^3$
for $\displaystyle \by_k = \sum_{j=1}^{k} {S}^{k-j}(S_j-S) {S}^{(j-1)}\bd_1$ and $\bx_0=(1,1,...,1)^T$ such that $S\bx_0=\bx_0$.
\end{itemize}
\end{theorem}
\proof
The proof follows the line of reasoning of the proof of \cite[Theorem 6]{DL95}.
For $\bd_1\in\RR^{N\times 3}$ we define
$$
\bu_{k+1,\ell}:=S^\ell S^{(k)} \bd_1,\quad \ell \geq 0,\quad  k\geq 0.
$$
From assumption $(i)$ we know that $\lim_{\ell\rightarrow +\infty} \bu_{k+1,\ell}$ is well defined. Next, with the notation $\bu_{k+1}:=\lim_{\ell\rightarrow +\infty} \bu_{k+1,\ell}$, we prove
 that the sequence $\{\bu_{k+1},\ k\geq 0\}$ is a Cauchy sequence. Indeed, in view of Proposition~\ref{prop:M_stable}, Proposition~\ref{lemmanew} and assumption $(iii)$ we have
$$
\|\bu_{k+1}-\bu_{k}\|_{\infty}=\left \|\lim_{\ell\rightarrow +\infty} S^\ell (S_{k}-S) S^{(k-1)} \bd_1 \right \|_{\infty} \leq \mathpzc{\bar C}\|S_{k}-S\|_{\infty} \|\bd_1\|_{\infty} \leq \frac{\mathpzc{\tilde C}}{\sigma^{k}},
$$
and thus, for $s\geq 1$,
\begin{equation}\label{eq:convrate_uk}
\|\bu_{k+s}-\bu_{k}\|_{\infty} \leq \sum_{j=1}^{s} \|\bu_{k+j}-\bu_{k+j-1}\|_{\infty} \leq \frac{\mathpzc{\hat C}}{\sigma^{k}} \sum_{j=0}^{s-1} \frac{1}{\sigma^{j}},
\end{equation}
with $\mathpzc{\bar C}$, $\mathpzc{\tilde C}$  and $\mathpzc{\hat C}$ finite positive constants.
Hence, the vector $\bu:=\lim_{k\rightarrow +\infty} \bu_{k}$ is well defined. \\
The next step is to show that  $\bu$ is in fact the limit of the sequence $\bd_m$, that is, $\bu=\lim_{m\rightarrow +\infty} S^{(m)} \bd_1$.
To this purpose, with the notation $\bd_{k+\ell+1}:=S^{(k+\ell)} \bd_1=S_{k+\ell}\cdots S_1\bd_1$,
we show that the difference $\| \bd_{k+\ell+1}-\bu_{k+1,\ell}\|_{\infty}$ becomes arbitrarily small for large enough $\ell$ and $k$.
Indeed,
$$
\| \bd_{k+\ell+1}-\bu_{k+1,\ell}\|_{\infty}=\left \|\left(S_{k+\ell}\cdots S_{k+1}-S^\ell\right)S_k\cdots S_1\bd_1 \right \|_{\infty} \leq
\mathpzc{C}\| S_{k+\ell}\cdots S_{k+1}-S^\ell\|_{\infty},
$$
with $\mathpzc{C}$ a finite positive constant.
In view of Proposition \ref{prop:Mk_ricorrenza} (with $S_{k+1}$ playing the role of $M_1$), using again $(iii)$ we arrive at
\begin{equation}\label{eq:convrate1}
\| \bd_{k+\ell+1}-\bu_{k+1,\ell}\|_{\infty} \leq  \mathpzc{\tilde C} \sum_{j=k+1}^{k+\ell} \| {S}^{k+\ell-j}(S_j-S) {S}^{(j-1)}\bd_1\|_{\infty} \leq \mathpzc{\hat C} \sum_{j=k+1}^{k+\ell} \frac{1}{\sigma^{j}},
\end{equation}
where again $\mathpzc{\tilde C}$, $\mathpzc{\hat C}$ are finite positive constants. In conclusion, for large enough $\ell$ and $k$, $\| \bd_{k+\ell+1}-\bu_{k+1,\ell}\|_{\infty} $ also becomes arbitrarily small since so does the right hand side of \eqref{eq:convrate1}.\\
Using the triangular inequality  we write
$$
\| \bd_{k+\ell+1}-\bu\|_{\infty} \le  \| \bd_{k+\ell+1}-\bu_{k+1,\ell}\|_{\infty} +  \| \bu_{k+1,\ell}-\bu\|_{\infty}
$$
and  conclude that, for large enough $\ell$ and $k$, $\|S^{(k+\ell)}\bd_1-\bu\|_{\infty}$ can be made arbitrarily small. In other words,
\begin{equation}\label{eq:u_lim}
\bu=\lim_{m\rightarrow +\infty} S^{(m)} \bd_1.
\end{equation}
We continue by showing that the vector $\bu$ is an eigenvector of $S$ associated with the eigenvalue $1$ (\ie, $S\bu=\bu$). Indeed, observing that  $S\bu_{k+1,\ell}=\bu_{k+1,\ell+1}$ we write
$$
\|  S\bu-\bu \|_{\infty} \leq  \|  S\bu-S\bu_{k+1,\ell} \|_{\infty} + \|  \bu_{k+1,\ell+1}- \bd_{k+\ell+2} \|_{\infty} + \| \bd_{k+\ell+2} -\bu\|_{\infty},
$$
with the right hand side that tends to $0$ for $k$ and $\ell$ going to $+\infty$.
In view of assumption $(i)$ and \eqref{eq:u_lim}, we can thus conclude convergence of the sequence
$$
\{\by_k, \ k\geq 0\}, \qquad \hbox{with} \qquad  \by_k:=S^{(k)}\bd_1-S^k\bd_1=\bd_{k+1}-S^k\bd_1=\sum_{j=1}^{k} {S}^{k-j}(S_j-S) {S}^{(j-1)}\bd_1.
$$
Moreover, denoting $\by:=\lim_{k \rightarrow  +\infty} \by_k$, from the fact that $S\bu=\bu$ we can also conclude that $S\by=\by$, which means that $\by$ lies in the eigenspace corresponding to the right eigenvector of $S$ associated to the eigenvalue $\lambda_0=1$.
Therefore it follows that $\by$ is of the form $\by=\bx_0\bbeta_0^T$ with $\bx_0=(1,1,...,1)^T$, which implies that $\bbeta_0$ can be written as
$\bbeta_0=\by^T \, \frac{\bx_0}{\bx_0^T \bx_0}=\by^T \, \frac{\bx_0}{N}$. \\
From \eqref{eq:global_cpset} we then write
\begin{equation}\label{eq:step1_conv}
\lim_{k\rightarrow +\infty}\bd_{k+1}=\lim_{k\rightarrow +\infty} {S}^k\bd_1+\bx_0\bbeta_0^T,
\end{equation}
and in view of Proposition \ref{prop:conv_stat}, after replacing \eqref{eq:lim} in equation  \eqref{eq:step1_conv}, we arrive at
\begin{equation}\label{eq:step1b_conv}
\lim_{k \rightarrow + \infty} \bd_{k+1} = \bx_0(\bq_0+\bbeta_0)^T, \quad \hbox{with} \quad \bq_0=\bd_1^T \tilde{\bx}_0.
\end{equation}
Then, taking into consideration assumption ({\it ii}) and Theorem \ref{prop:asympt_equiv}, we have that
\begin{equation}\label{eq:unif_conv_phi}
\lim_{k \rightarrow + \infty} \sup_{(u,v) \in \Omega_{k+1}} \|\bPhi_{k+1}(u,v)-\overline{\bPhi}(u,v)\|_{\infty}=0.
\end{equation}
The above means that, for any $\varepsilon>0$ and for large enough $k$,
\begin{equation}\label{eq_uniformboundPhi}
\sup_{(u,v) \in \Omega_{k+1}} \|\bPhi_{k+1}(u,v)\|_\infty \leq \sup_{(u,v) \in \Omega_{k+1}} \|\overline{\bPhi}(u,v)\|_{\infty}+ \varepsilon \leq \sup_{(u,v) \in \Omega} \|\overline{\bPhi}(u,v)\|_{\infty} + \varepsilon,
\end{equation}
i.e., $\sup_{(u,v) \in \Omega_{k+1}} \|\bPhi_{k+1}(u,v)\|_\infty$ is uniformly bounded.
After recalling that $\overline{\bPhi}(u,v)^T \bx_0=1$ for all $(u,v) \in \Omega_{k+1}$ (in light of the arguments in Remark \ref{rem:POU}), we continue by writing, for all $j \in \mathbb{J}_{3n}$,
$$
\hspace{-0.3cm}
\begin{array}{l}
\displaystyle
\sup_{(u,v) \in \omega_{k+1}^{[j]}} \hspace{-0.1cm} \| \br_{k+1}^{[j]}(u,v)^T-(\bq_0+\bbeta_0)^T\|_{\infty}=\\
\displaystyle
\sup_{(u,v) \in \omega_{k+1}^{[j]}} \hspace{-0.1cm} \| \br_{k+1}^{[j]}(u,v)^T-\overline{\bPhi}(u,v)^T \bx_0(\bq_0+\bbeta_0)^T\|_{\infty}=\\
\displaystyle
\sup_{(u,v) \in \omega_{k+1}^{[j]}} \hspace{-0.1cm} \| \bPhi_{k+1}(u,v)^T \bd_{k+1}-\overline{\bPhi}(u,v)^T \bx_0(\bq_0+\bbeta_0)^T\|_{\infty}=\\
\displaystyle
\sup_{(u,v) \in \omega_{k+1}^{[j]}} \hspace{-0.1cm} \| \bPhi_{k+1}(u,v)^T \bd_{k+1}-\bPhi_{k+1}(u,v)^T \bx_0(\bq_0+\bbeta_0)^T + \\ \qquad \qquad \bPhi_{k+1}(u,v)^T \bx_0(\bq_0+\bbeta_0)^T    -\overline{\bPhi}(u,v)^T \bx_0(\bq_0+\bbeta_0)^T\|_{\infty} \leq \smallskip \smallskip \\
\\
\displaystyle
\sup_{(u,v) \in \omega_{k+1}^{[j]}} \hspace{-0.1cm} \| \bPhi_{k+1}(u,v)^T\|_{\infty} \, \|\bd_{k+1}- \bx_0(\bq_0+\bbeta_0)^T \|_{\infty} +\\
\displaystyle \sup_{(u,v) \in \omega_{k+1}^{[j]}} \hspace{-0.1cm} \| \bPhi_{k+1}(u,v)^T-\overline{\bPhi}(u,v)^T \|_{\infty} \, \|\bx_0(\bq_0+\bbeta_0)^T\|_{\infty}.
\end{array}
$$
Since $\displaystyle \lim_{k \rightarrow +\infty} \|\bd_{k+1}- \bx_0(\bq_0+\bbeta_0)^T \|_{\infty}=0$ (in light of \eqref{eq:step1b_conv}),
$\displaystyle \lim_{k \rightarrow +\infty} \sup_{(u,v) \in \Omega_{k+1}} \| \bPhi_{k+1}(u,v)^T-\overline{\bPhi}(u,v)^T \|_{\infty}=0$ (in light of \eqref{eq:unif_conv_phi}), $\displaystyle \sup_{(u,v) \in \Omega_{k+1}} \| \bPhi_{k+1}(u,v)^T\|_{\infty}$ is uniformly bounded (in light of \eqref{eq_uniformboundPhi}) and $\|\bx_0(\bq_0+\bbeta_0)^T\|_{\infty}$ is bounded, we finally obtain
$$
\lim_{k \rightarrow +\infty} \sup_{(u,v) \in \omega_{k+1}^{[j]}} \| \br_{k+1}^{[j]}(u,v)^T-(\bq_0+\bbeta_0)^T\|_{\infty}=0, \quad \forall \ j \in \mathbb{J}_{3n},
$$
which concludes the proof.
\endproof

\smallskip
\begin{remark}
It is worthwhile to stress that, by requiring that the matrix sequence $\{S_k, k \geq 1 \}$ (identifying a non-stationary subdivision scheme in the vicinity of an extraordinary element) converges towards the matrix $S$ of a convergent stationary scheme  faster than $\frac{1}{\sigma^k}$,  $\sigma>1$, the convergence of the non-stationary subdivision scheme follows.
\end{remark}

\smallskip
Now, following the notation in \cite{PR08}, we denote with  $\lambda_r$, $r=0, \ldots, \overline{r}$, $0 \leq \overline{r}\leq N-1$, the $\overline{r}+1$ different eigenvalues of $S \in \RR^{N \times N}$ sorted in decreasing order according to their magnitude, i.e.,  $|\lambda_0|\geq |\lambda_1| \geq \ldots \geq |\lambda_{\overline{r}}|$. Moreover, for $r \geq 0$ we denote by $\ell_r+1$ the algebraic multiplicity of $\lambda_r$.
As emphasized in Remark \ref{rem:convergenceimplies}, it is a known fact that, for a rotationally symmetric, convergent stationary scheme $\bar{\mathscr S}$  associated with $S$, all $\overline{r}+1$ eigenvalues have magnitude less than 1, except $\lambda_0$ which is required to be exactly $1$ and with algebraic and geometric multiplicity $1$. It means that
$1=\lambda_0 > |\lambda_1| \geq \ldots \geq |\lambda_{\overline{r}}|$ and $\ell_0=0$.
Moreover, the eigenvector associated to the unique dominant eigenvalue $\lambda_0=1$ is required to be $\bx_0=(1,1,...,1)^T\in \RR^{N}$ (see, e.g., \cite{PR08,R95,Z00}).
Thus, exploiting the Jordan decomposition of $S^k$ and the equality $S^k \bd_1=XJ^k\bX^{-1} \bd_1$, we can write
\begin{equation}\label{eq:convrate_stat}
S^{k} \bd_1= \bx_0 \bq_0^T + O(|\lambda_1|^k {\bf 1}),
\end{equation}
where, with a slight abuse of notation, $O(|\lambda_1|^k {\bf 1})$ denotes a vector in $\RR^{N\times 3}$ with all its entries behaving as $O(|\lambda_1|^k)$.

\smallskip \noindent
Equation \eqref{eq:convrate_stat} implies the following convergence rate result for the sequence $\{\by_k, \ k \geq 0\}$.

\begin{corollary}\label{convergence_rateNEW}
Let  $\bar{\mathscr S}$  be a symmetric, convergent, stationary scheme represented by $S$ and denote by $\lambda_1$ the subdominant eigenvalue of $S$.
Under the assumptions of Theorem \ref{theo:convergenceNEW} with the additional requirement that $\sigma>\frac{1}{|\lambda_1|}>1$, for
$
\bu=\lim_{k\rightarrow +\infty}S^{(k)}\bd_1
$
we have
\begin{equation}\label{eq:convrate_u}
S^{(k)}\bd_1=\bu+ O\left( \frac{1}{\sigma^{k}} {\bf 1}\right).
\end{equation}
Consequently, $\by_k$ as in \eqref{eq:global_cpset} satisfies
\begin{equation}\label{eq:convrate}
\by_k=\bx_0 \bbeta_0^T + O\left(\frac{1}{\sigma^k} {\bf 1}\right).
\end{equation}
\end{corollary}
\proof
We use the notation of the proof of Theorem \ref{theo:convergenceNEW} and start proving \eqref{eq:convrate_u}.
First we write
$$\bu-S^{(k+\ell)} \bd_1=(\bu-\bu_{k+1})+(\bu_{k+1}-\bu_{k+1,\ell})+(\bu_{k+1,\ell}-S^{(k+\ell)} \bd_1).$$
Then, by \eqref{eq:convrate_uk}, \eqref{eq:convrate1} and \eqref{eq:convrate_stat}  we obtain
\begin{align*}
\|\bu-S^{(k+\ell)} \bd_1\|_{\infty} & \leq \|\bu-\bu_{k+1}\|_{\infty}+
\|\bu_{k+1}-\bu_{k+1,\ell}\|_{\infty}+
\|\bu_{k+1,\ell}-S^{(k+\ell)} \bd_1\|_{\infty} \\
& \leq \mathpzc{C}_1 |\lambda_1|^{\ell} + \frac{\mathpzc{C}_2}{{\sigma}^{k}},
\end{align*}
with $\mathpzc{C}_1$, $\mathpzc{C}_2$ finite positive constants.
Since $|\lambda_1|<1$, we can find $\bar{L}$ such that $|\lambda_1|^{\ell}<\frac{1}{\sigma^k}$ for all $\ell>\bar{L}$. Therefore,
$$
\|\bu-S^{(k+\bar{L}+1)} \bd_1\|_{\infty} \leq \frac{\mathpzc{C}_3}{{\sigma}^{k}},
$$
with $\mathpzc{C}_3$ a finite positive constant. Hence, taking the limit to $+ \infty$ with respect to $k$, \eqref{eq:convrate_u} follows.\\
Similarly, since $\bx_0 (\bq_0 + \bbeta_0)^T=\lim_{k \rightarrow +\infty} \bd_k=\bu$ due to \eqref{eq:u_lim}, to prove the result in \eqref{eq:convrate} we write
$$
\by_{k+\ell}-\bx_0 \bbeta_0^T=(S^{(k+\ell)} \bd_1-\bu)-(S^{k+\ell}\bd_1-\bx_0 \bq_0^T),
$$
and consider the triangular inequality
\begin{align*}
\|\by_{k+\ell}-\bx_0 \bbeta_0^T\|_{\infty} & \leq \|S^{(k+\ell)} \bd_1-\bu_{k+1,\ell}\|_{\infty}+\|\bu_{k+1,\ell}-\bu_{k+1}\|_{\infty}\\
&+ \|\bu_{k+1}-\bu\|_{\infty} +\|S^{k+\ell}\bd_1-\bx_0 \bq_0^T\|_{\infty}.
\end{align*}
Using again \eqref{eq:convrate_uk}, \eqref{eq:convrate1} and \eqref{eq:convrate_stat}, we obtain the upper bound
$$
\|\by_{k+\ell}-\bx_0 \bbeta_0^T\|_{\infty} \leq \tilde{\mathpzc{C}}_1 |\lambda_1|^{\ell} + \frac{\tilde{\mathpzc{C}}_2}{{\sigma}^{k}}+\tilde{\mathpzc{C}}_3 |\lambda_1|^{k+\ell}.
$$
Then, applying the same reasoning as before, \eqref{eq:convrate} is proven.
\endproof

\subsection{Normal continuity analysis at the limit point of an extraordinary element}\label{sec:4}
Aim of this section is to provide sufficient conditions to show that a rotationally symmetric, convergent, non-stationary subdivision scheme ${\mathscr S}$ produces a normal continuous surface at the limit point of an extraordinary element.\\

\noindent
For the rotationally symmetric, stationary subdivision scheme $\bar{\mathscr S}$ we assume all ineffective eigenvectors of the associated local subdivision matrix $S$ to be ruled out
(according to \cite[Theorem 4.20]{PR08} the absence of ineffective eigenvectors can
be required without loss of generality) and
the ordered eigenvalues of $S$ to satisfy
$$
1=\lambda_0> \lambda_1 > |\lambda_2| \quad \hbox{with} \quad \lambda_1 \in \RR^+, \, \ell_1=1,
$$
i.e., the sub-dominant eigenvalue $\lambda_1$ to be real, double and with geometric multiplicity equal to algebraic multiplicity (namely, $\lambda_1$ non-defective). In this case, the eigenvectors associated to $\lambda_1$ are linearly independent. In the following we denote by
$\bx_0=(1,1,...,1)^T\in \RR^{N}$ the eigenvector associated to $\lambda_0=1$, and by  $\bx_1^0, \bx_1^1 \in \RR^{N}$ the two linearly independent eigenvectors associated to $\lambda_1$.
Moreover, we assume that the stationary scheme $\bar{\mathscr S}$ is $C^1$-convergent in regular regions and
the planar ring defined by $\overline{\bPsi}(u,v)^T=\overline{\bPhi}(u,v)^T \left(\bx_1^0, \bx_1^1\right) \in \RR^{1\times 2}$ (where $\overline{\bPhi}(u,v)\in \RR^{N}$ denotes the associated basic limit function vector) is such that
\begin{equation}\label{eq:cond_charmap}
sign \left ( \det(\cJ \overline{\bPsi}(u,v)^T) \right)  \; \hbox{is non-zero and constant for all} \; (u,v) \in \Omega_1,
\end{equation}
with
$$
\cJ \overline{\bPsi}(u,v)^T:=\begin{pmatrix}\partial_u \overline{\bPsi}(u,v)^T\\\partial_v \overline{\bPsi}(u,v)^T\end{pmatrix}\in \RR^{2\times 2}.$$

\begin{remark}
Note that stationary subdivision schemes that possess a real, double subdominant
eigenvalue are those commonly termed \emph{standard} (see, e.g., \cite[Chapter 5.2]{PR08}). We restrict our attention to them since they are the ones of practical relevance. Moreover, note that the assumption in (\ref{eq:cond_charmap}) is nothing but the notion of regularity of the characteristic map of $\bar{\mathscr S}$ (see, e.g., \cite{PR08} for details).
\end{remark}
In the following we provide sufficient conditions to show that a rotationally symmetric, non-stationary subdivision scheme ${\mathscr S}$ produces a normal continuous surface at the limit point $\br_c=\bq_0+\bbeta_0$ (see Definition~\ref{def:G1}).

\begin{theorem}\label{teo:G1}
Let ${\mathscr S}$ be a rotationally symmetric, non-singular, non-stationary subdivision scheme whose action in an irregular region is described by a matrix sequence  $\{S_k, k \geq 1\}$ and whose action in the regular patch ring $\br_{k+1}$ is described by a basic limit function vector $\bPhi_{k+1}(u,v)$ that satisfies the condition $\bPhi_{k+1}(u,v)^T {\bf 1}=\alpha_{k+1} \in \RR$.
Moreover, let $\bar{\mathscr S}$ be a rotationally symmetric, standard, stationary subdivision scheme that in the same irregular region is associated with $S$.
Assume that:

\begin{itemize}
\item[(i)] $\bar{\mathscr S}$ is $C^1$-convergent in regular regions, with a symbol $c(\bz)$ containing the factor $(1+z_1)(1+z_2)$,  and satisfies \eqref{eq:cond_charmap};
\item[(ii)] in regular regions ${\mathscr S}$ is defined by the symbols $\{c^{(k)}(\bz), \, k \geq 1\}$ where each  $c^{(k)}(\bz)$ contains the factor $(1+z_1)(1+z_2)$;
    \item[(iii)] in regular regions ${\mathscr S}$ is asymptotically equivalent of order $1$ to $\bar{\mathscr S}$;
\item[(iv)] in the irregular region the matrices $S_k$ and $S$ satisfy, for all $k \geq 1$, $\|S_k-S\|_{\infty} \leq \frac{\mathpzc{C}}{\sigma^k}$
    with $\mathpzc{C}$ some finite positive constant, $\sigma>\frac{1}{\lambda_1}>1$.
\end{itemize}
Then the subdivision surface generated by ${\mathscr S}$ is normal continuous at the limit point $\br_c$.
\end{theorem}

\proof
First we observe that from $(i)$  and Remark  \ref{rem:convergenceimplies} the matrix $S$ has a simple dominant eigenvalue $\lambda_0=1$. Also, from Theorem \ref{prop:asympt_equiv_order1} we know that ${\mathscr S}$ is $C^1$-convergent in regular regions and from Theorem \ref{theo:convergenceNEW} we also know that ${\mathscr S}$ is convergent in the irregular region.
Therefore, the two sequences $\{\br_{k+1}, k \geq 0\}$ and $\{\bd_{k+1}, k \geq 0\}$ converge.
To simplify the analysis, we do not consider the full expression of a sequence of rings, but only the asymptotic behavior of the dominant terms as $k$ tends to infinity.
Due to assumption ({\it i}), the eigenvalues of $S$ satisfy $1=\lambda_0> \lambda_1 > |\lambda_i| ,\ i=2, \ldots, \overline{r}$ and the sub-dominant eigenvalue $\lambda_1$ has geometric multiplicity and algebraic multiplicity two \cite{PR08}.
Thus, recalling from Proposition \ref{prop:Mk_ricorrenza}  that
$$\bd_{k+1}=S^{k} \bd_1+ \by_k \qquad \hbox{with} \qquad \by_k=\sum_{j=1}^{k}S^{k-j}(S_j-S)S^{(j-1)} \bd_1$$
and exploiting the Jordan decomposition of $S^k$ given by $S^k=XJ^k\bX^{-1}$,
we obtain
$$
S^{k} \bd_1= \displaystyle \bx_0 \bq_0^T + \lambda_1^k (\bx_1^0 (\bq_1^0)^T+\bx_1^1 (\bq_1^1)^T) + o\left(\lambda_1^{k} \mathbf{1}\right),
$$
with $\bx_1^0$ and $\bx_1^1$ denoting the two linearly independent eigenvectors associated to $\lambda_1$, $\bq_1^0$, $\bq_1^1$ two vectors in $\RR^3$ and $o\left(\lambda_1^{k} \mathbf{1}\right)$ a vector in $\RR^{N\times 3}$ with all its entries behaving as $o\left(\lambda_1^{k}\right)$.
Since $\frac{1}{\sigma}<\lambda_1$, in view of Corollary \ref{convergence_rateNEW} we also have  that
$$
\by_k=\bx_0 \bbeta_0^T + o\left(\lambda_1^{k} \mathbf{1}\right).
$$
This yields
\begin{equation}\label{eq:d}
\displaystyle \bd_{k+1}= \displaystyle \bx_0 (\bq_0^T+\bbeta_0^T)+ \lambda_1^k (\bx_1^0 (\bq_1^0)^T+\bx_1^1 (\bq_1^1)^T)
+ o\left(\lambda_1^{k} \mathbf{1}\right).
\end{equation}
\noindent
Parameterizing the regular patch ring $\br_{k+1}$ using the basic limit function vector $\bPhi_{k+1}$, we can write $(\br_{k+1}^{[j]})^T$, for each $j \in \mathbb{J}_{3n}$, as (cf. Equation \eqref{eq:r_k})
$$
(\br_{k+1}^{[j]}(u,v))^T=\bPhi_{k+1}^T(u,v) \bd_{k+1},\quad (u,v)\in \omega^{[j]}_{k+1},\quad j \in \mathbb{J}_{3n}.
$$
Using Remark \ref{rem:POU} and introducing the shorthand notation $\alpha_{k+1}$ for the value $\bPhi_{k+1}(u,v)^T \bx_0 \in \RR$, thanks to \eqref{eq:d}, we have
\begin{align*}
(\br_{k+1}^{[j]}(u,v))^T&=\bPhi_{k+1}(u,v)^T \left ( \bx_0 (\bq_0^T+\bbeta_0^T)+ \lambda_1^k  (\bx_1^0 (\bq_1^0)^T+\bx_1^1 (\bq_1^1)^T)
+  o\left(\lambda_1^{k} \mathbf{1}\right)
\right)\\
&=\alpha_{k+1} (\bq_0^T+\bbeta_0^T)+ \lambda_1^k \bPhi_{k+1}(u,v)^T (\bx_1^0 (\bq_1^0)^T+\bx_1^1 (\bq_1^1)^T)\\
&+ \bPhi_{k+1}(u,v)^T o\left(\lambda_1^{k} \mathbf{1}\right).
\end{align*}
\noindent
To verify the normal continuity of the limit surface at the limit point $\br_c=\bq_0+\bbeta_0$, we first observe, using Remark \ref{rem:POU}, that
$$
\partial_u \alpha_{k+1}=\partial_v \alpha_{k+1}=0,
$$
and then write
\begin{align*}
\partial_u (\br_{k+1}^{[j]}(u,v))^T& =\partial_u \Big( \alpha_{k+1} (\bq_0^T+\bbeta_0^T)+ \lambda_1^k \bPhi_{k+1}(u,v)^T (\bx_1^0 (\bq_1^0)^T+\bx_1^1 (\bq_1^1)^T)\\
& \qquad + \bPhi_{k+1}(u,v)^T  o\left(\lambda_1^{k} \mathbf{1}\right) \Big)\\
&=\lambda_1^k \partial_u  \bPhi_{k+1}(u,v)^T (\bx_1^0 (\bq_1^0)^T+\bx_1^1 (\bq_1^1)^T) + \partial_u  \bPhi_{k+1}(u,v)^T
 o\left(\lambda_1^{k} \mathbf{1}\right)\\
&=\lambda_1^k \partial_u  \bPhi_{k+1}(u,v)^T \left(  (\bx_1^0 (\bq_1^0)^T+\bx_1^1 (\bq_1^1)^T) + \frac{ o\left(\lambda_1^{k} \mathbf{1}\right)}{\lambda_1^k}  \right)
\end{align*}
and, similarly,
$$
\partial_v (\br_{k+1}^{[j]}(u,v))^T=\lambda_1^k \partial_v  \bPhi_{k+1}(u,v)^T \left( (\bx_1^0 (\bq_1^0)^T+\bx_1^1 (\bq_1^1)^T) +  \frac{ o\left(\lambda_1^{k} \mathbf{1}\right) }{\lambda_1^k}   \right).
$$
Since the $\bPhi_{k+1}$ and their derivatives converge uniformly to $\overline{\bPhi}$ and its derivatives due to Theorem \ref{prop:asympt_equiv} and Corollary \ref{cor:unif_conv_der}, we see that the $\lambda_1^{-k} \br_{k}$ and their derivatives converge uniformly to the $C^1$ function
$$\overline{\bPhi}^T \left(\bx_1^0 (\bq_1^0)^T+\bx_1^1 (\bq_1^1)^T \right)$$
which maps into the linear space spanned by $\bq_1^0$ and $\bq_1^1$.
Consequently the normal vectors of the $\br_{k}$, which are the same as of the scaled functions $\lambda_1^{-k} \br_{k}$, converge uniformly to
$$
\bn_{\infty}:=\frac{(\bq_1^0)^T \wedge (\bq_1^1)^T}{\| (\bq_1^0)^T \wedge (\bq_1^1)^T  \|_2}.
$$
The latter shows that the limit surface $\br$ obtained by the non-stationary subdivision scheme ${\mathscr S}$ is
normal continuous at the limit point $\br_c$, which concludes the proof.
\endproof

\medskip
\begin{remark}
Theorem \ref{teo:G1} provides a sufficient condition for verifying the normal continuity of a non-stationary scheme.
It consists in verifying that the matrix sequence $\{S_k, k \geq 1 \}$ (identifying the non-stationary scheme
in the vicinity of an extraordinary element) converges towards a $C^1$-regular, standard, stationary scheme (identified with the matrix $S$) faster than $\lambda_1^k$ (i.e., with convergence rate $o(\lambda_1^k)$), where $\lambda_1$ denotes the real, double subdominant eigenvalue of $S$.
\end{remark}

\section{Application examples}\label{sec:5}
In this section we use Theorem \ref{theo:convergenceNEW} to study the convergence of two non-stationary subdivision schemes, defined on quadrilateral meshes, in the neighborhood of extraordinary elements. Also,  we use  Theorem \ref{teo:G1} to prove that the limit surfaces obtained by such schemes are normal continuous at the limit points of the corresponding extraordinary elements. This partially proves a conjecture given in  \cite[Section 5]{FMW14} where only numerical evidence for $C^1$-regularity was shown.

\subsection{Generalized trigonometric spline surfaces of order $3$}

\begin{figure}[h!]
\centering
\includegraphics[trim={0.15cm 0.05cm 0.25cm 0.05cm}, clip, height=0.33\textwidth]{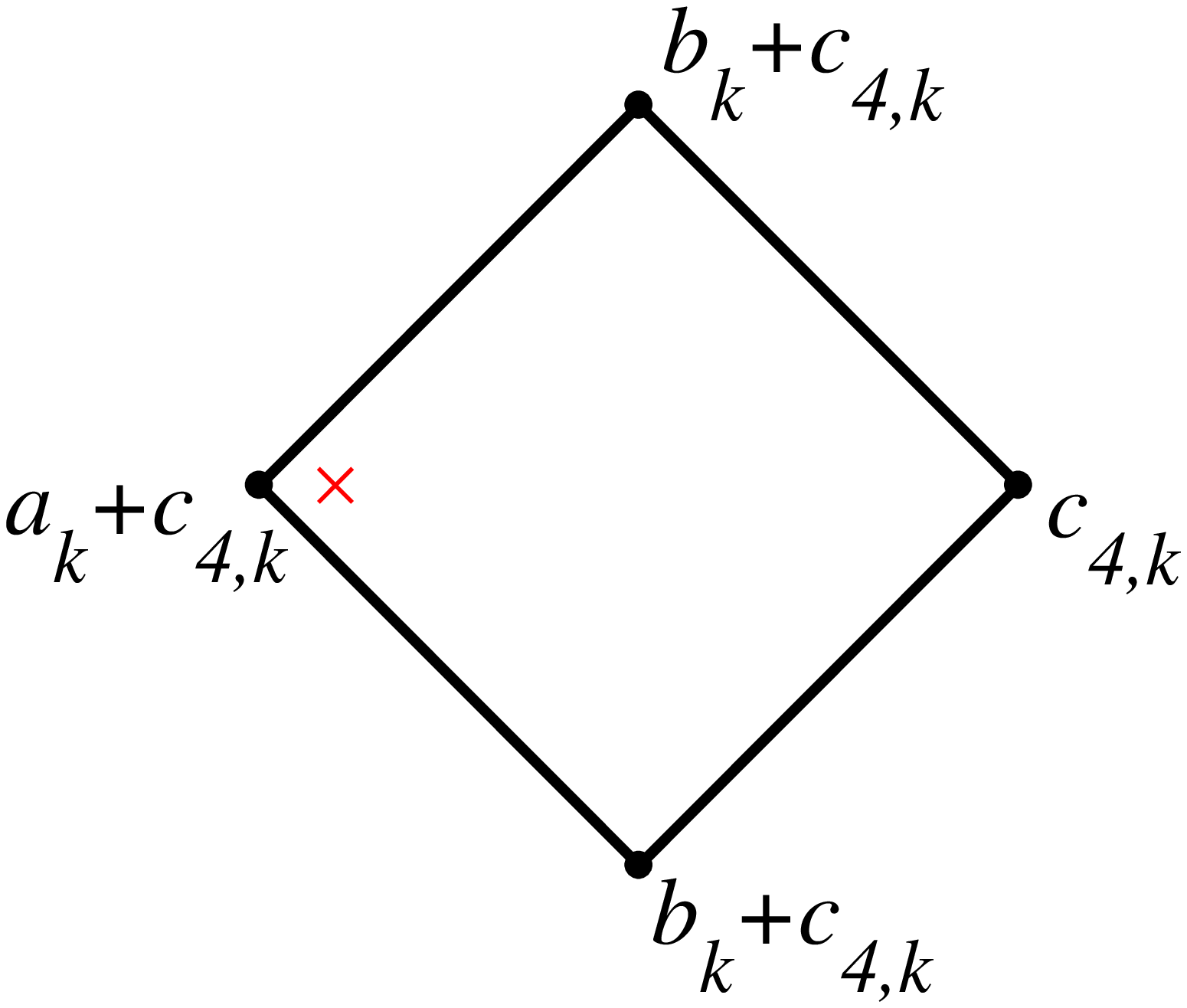}\hspace{0.85cm}
\includegraphics[trim={0.25cm 0.05cm 0.25cm 0.05cm}, clip, height=0.3\textwidth]{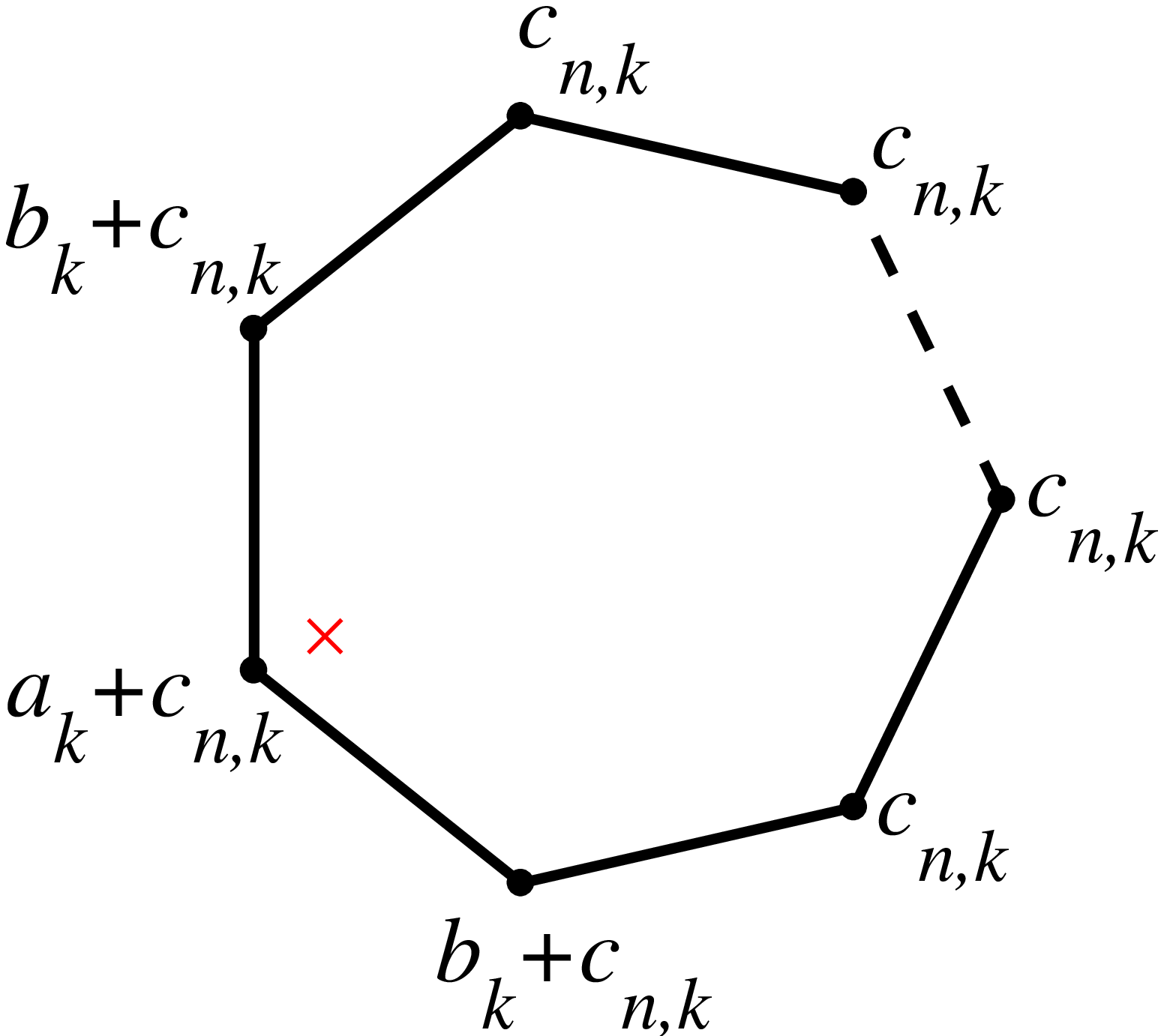}
\caption{Graphical illustration of the $k$-th level geometric refinement rules of the non-stationary subdivision scheme generalizing trigonometric spline surfaces of order $3$. The red cross represents the new point inserted by the geometric refinement rules in the case of regular (left) and extraordinary (right) faces. The weights appearing in the refinement rules are the ones specified in  \eqref{eq:coeff1} and \eqref{eq:coeff2}. (Color figure online.)}
\label{fig:ds_stencils}
\end{figure}

In \cite{JSD02}, the authors presented a non-stationary subdivision scheme which produces tensor-product  trigonometric spline surfaces of order 3 except in the neighborhood of extraordinary faces. This non-stationary scheme can be seen as a generalization of the well-known stationary Doo-Sabin scheme \cite{DS78} yielding polynomial spline surfaces of order 3 except in the neighborhood of extraordinary faces.
Figure \ref{fig:ds_stencils} illustrates the $k$-th level geometric refinement rules of this non-stationary scheme. We do not include a figure illustrating the topologic refinement rules since they are exactly the same as the ones used by the standard (stationary) Doo-Sabin scheme.

\begin{figure}[h!]
\centering
\subfigure[]{\includegraphics[trim={1.65cm 1.5cm 1.65cm 1.3cm}, clip, height=0.23\textwidth]{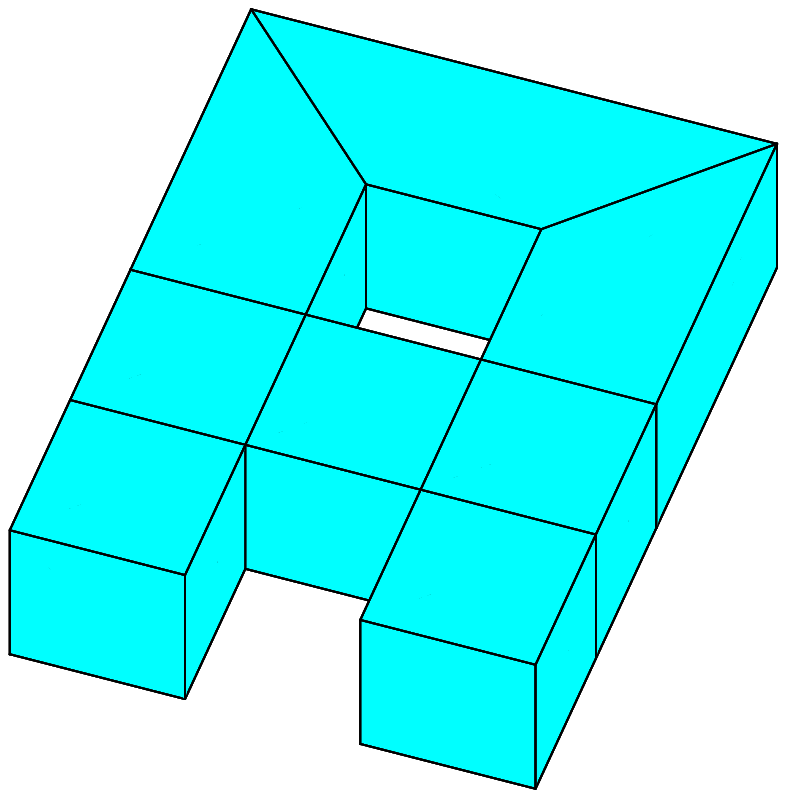}}
\subfigure[$h=\frac{1}{16}$]{\includegraphics[trim={2.1cm 1.8cm 2.1cm 2cm}, clip, height=0.23\textwidth]{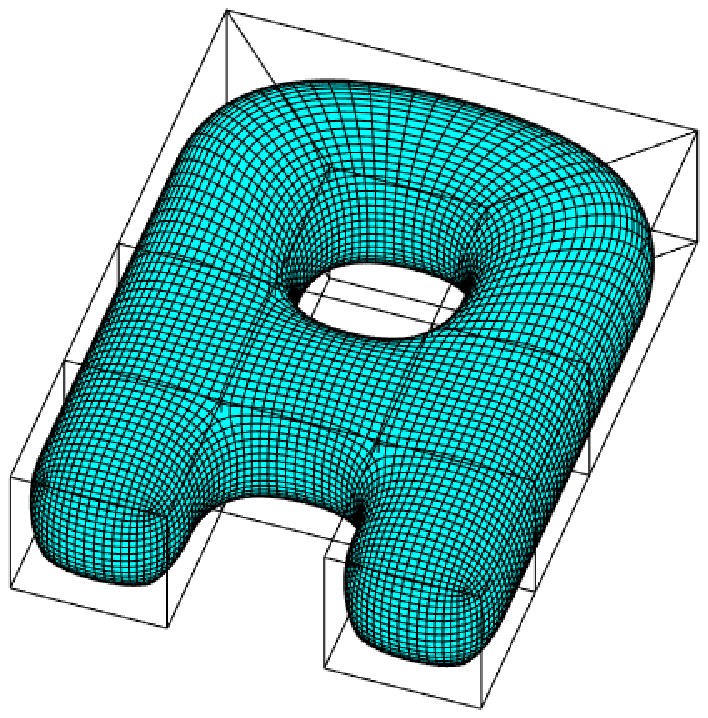}}
\subfigure[$h=1$]{\includegraphics[trim={2.1cm 1.8cm 2.1cm 2cm}, clip, height=0.23\textwidth]{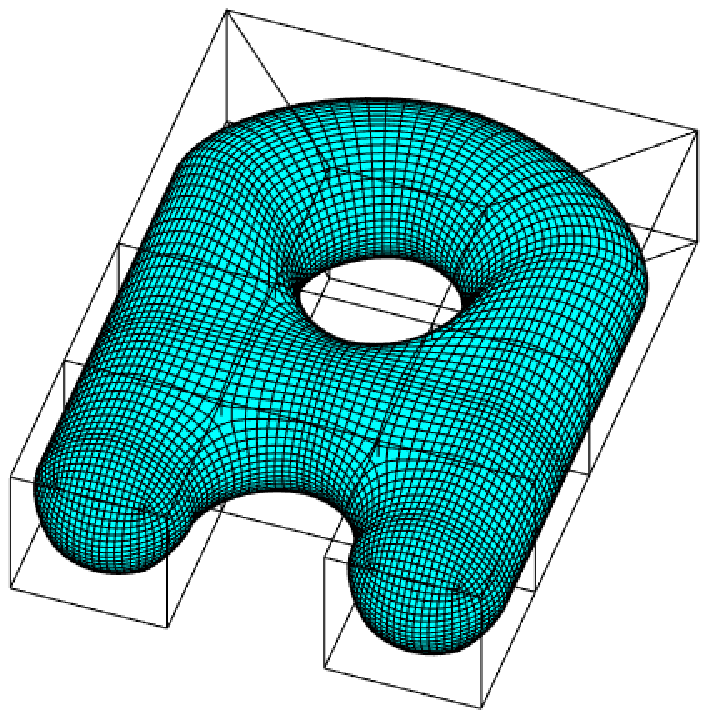}}
\vspace{0.3cm}
\setcounter{subfigure}{0}
\subfigure[]{\includegraphics[trim={2.4cm 2.1cm 2.4cm 2.1cm}, clip, height=0.22\textwidth]{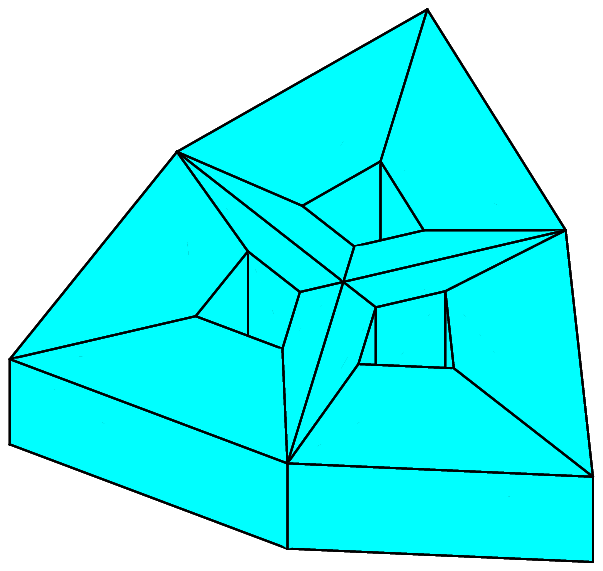}}
\subfigure[$h=\frac{1}{16}$]{\includegraphics[trim={2.8cm 2.5cm 2.8cm 2.5cm}, clip, height=0.22\textwidth]{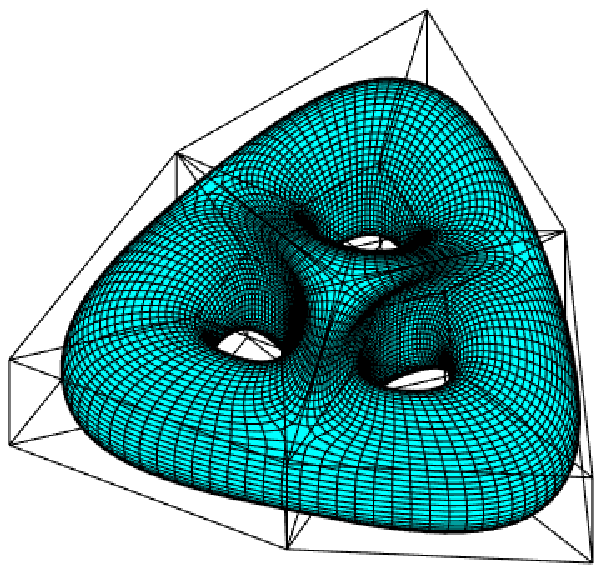}}
\subfigure[$h=1$]{\includegraphics[trim={2.8cm 2.5cm 2.8cm 2.5cm}, clip, height=0.22\textwidth]{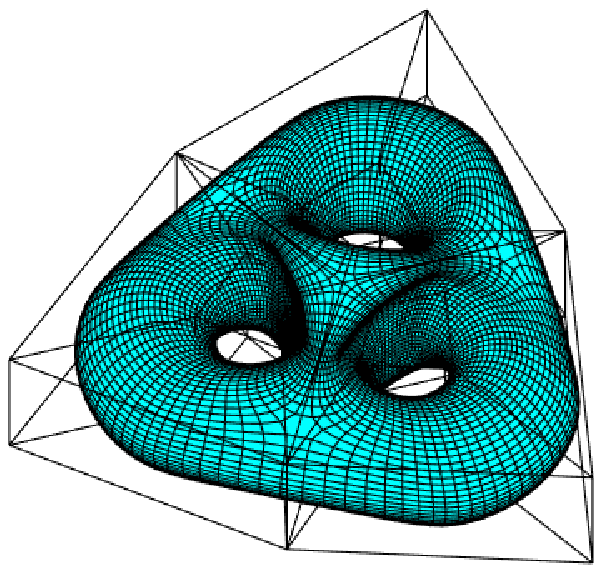}}
\caption{Original control mesh (a) and refined meshes (b,c) obtained by applying four iterations of the normalized non-stationary subdivision scheme generalizing trigonometric spline surfaces of order $3$ with two different values of $h \in [0, \frac{\pi}{3})$. (Color figure online.)}
\label{fig:ds_examples}
\end{figure}

\noindent
In regular regions, Doo-Sabin scheme is described by the subdivision mask
\begin{equation}\label{eq:mask_ds}
\bc=\begin{pmatrix}
\frac{1}{16} & \frac{3}{16} & \frac{3}{16} & \frac{1}{16} \smallskip \\
\frac{3}{16} & \frac{9}{16} & \frac{9}{16} & \frac{3}{16} \smallskip\\
\frac{3}{16} & \frac{9}{16} & \frac{9}{16} & \frac{3}{16} \smallskip\\
\frac{1}{16} & \frac{3}{16} & \frac{3}{16} & \frac{1}{16}
\end{pmatrix},
\end{equation}
while in irregular regions the refinement rules are written in terms of a subdivision matrix $S$ having the structure in  \eqref{eq:Sk_block_dual} with blocks
\begin{equation}\label{eq:ds}
\begin{array}{c}
B_0=\begin{pmatrix}
\frac{1}{4n}+\frac12 & 0 & 0 & 0 \smallskip\\
\frac{9}{16} & \frac{3}{16} & 0 & 0 \smallskip\\
\frac{9}{16} & \frac{3}{16} &\frac{1}{16} & \frac{3}{16} \smallskip\\
\frac{9}{16} & 0 & 0 & \frac{3}{16}
\end{pmatrix}, \qquad
B_1=\begin{pmatrix}
\frac{1}{4n}+\frac18 & 0 & 0 & 0 \smallskip\\
 0 & 0 & 0 & 0 \smallskip\\
0 & 0 & 0 & 0 \smallskip\\
\frac{3}{16} & \frac{1}{16} & 0 & 0
\end{pmatrix},
\medskip \\
B_i=\begin{pmatrix}
\frac{1}{4n} & 0 & 0 & 0 \smallskip\\
 0 & 0 & 0 & 0 \smallskip\\
0 & 0 & 0 & 0 \smallskip\\
0 & 0 & 0 & 0
\end{pmatrix}, \; \; i=2, \ldots, n-2, \qquad
B_{n-1}=\begin{pmatrix}
\frac{1}{4n}+\frac18 & 0 & 0 & 0 \smallskip\\
 \frac{3}{16} & 0 & 0 & \frac{1}{16} \smallskip\\
0 & 0 & 0 & 0 \smallskip\\
0 & 0 & 0 & 0
\end{pmatrix}.
\end{array}
\end{equation}
It is a well-known fact that Doo-Sabin scheme is convergent both in regular regions and in irregular regions, and the limit surface is $C^1$. Moreover, in regular regions the associated subdivision symbol is
$$c(z_1,z_2)=\frac{(z_1 + 1)^3(z_2 + 1)^3}{16},$$
which contains the factor $(1+z_1)(1+z_2)$. Thus it satisfies assumption ({\it i}) of Theorem \ref{theo:convergenceNEW} and assumption ({\it i}) of Theorem \ref{teo:G1}.\\
In regular regions, the non-stationary scheme in \cite{JSD02} is described  by the $k$-th level mask
\begin{equation}\label{eq:mask_ds_ns}
\bc^{(k)}=\begin{pmatrix}
c_{4,k} & b_{k}+c_{4,k} &  b_{k}+c_{4,k} & c_{4,k} \smallskip\\
b_{k}+c_{4,k} & a_{k}+c_{4,k} &  a_{k}+c_{4,k} &  b_{k}+c_{4,k} \smallskip\\
b_{k}+c_{4,k} & a_{k}+c_{4,k} &  a_{k}+c_{4,k} &  b_{k}+c_{4,k} \smallskip\\
c_{4,k} &  b_{k}+c_{4,k} &  b_{k}+c_{4,k} & c_{4,k}
\end{pmatrix}, \qquad k \geq 1,
\end{equation}
where for $h \in \left[0,\frac{\pi}{3}\right)$,
\begin{equation}\label{eq:coeff1}
c_{n,k}=\frac{1}{4n \cos^2\left(\frac{h}{2^k}\right)\cos^2\left(\frac{h}{2^{k-1}}\right)}, \quad n\in \NN, \ n \geq 4, \quad k\geq 1,
\end{equation}
and
\begin{equation}\label{eq:coeff2}
a_{k}=\frac{1}{4\cos^2\left(\frac{h}{2^k}\right)\cos\left(\frac{h}{2^{k-1}}\right)}+\frac{1}{4\cos^2\left(\frac{h}{2^k}\right)}, \ \  b_{k}=\frac{1}{8\cos^2\left(\frac{h}{2^k}\right)\cos\left(\frac{h}{2^{k-1}}\right)}, \quad k\geq 1.
\end{equation}
Therefore the associated subdivision symbol is
$$
c^{(k)}(z_1,z_2)=
\resizebox{10.5cm}{!}{$\frac{e^{{\rm i} \frac{h}{2^{k-1}}} (z_1 + 1) (z_2 + 1) (z_1 + e^{{\rm i} \frac{h}{2^{k-1}}}) (z_1 e^{{\rm i} \frac{h}{2^{k-1}}} + 1) (z_2 + e^{{\rm i} \frac{h}{2^{k-1}}}) (z_2 e^{{\rm i} \frac{h}{2^{k-1}}} + 1)}{(e^{{\rm i} \frac{h}{2^{k-2}}}+1)^2 (e^{{\rm i} \frac{h}{2^{k-1}}}+1)^2}
$},
$$
which contains the factor $(1+z_1)(1+z_2)$, thus satisfying assumption ({\it ii}) of Theorem \ref{teo:G1}.
Differently, in irregular regions the refinement rules are given in terms of the $k$-th level matrix $S_k$ having the structure in  \eqref{eq:Sk_block_dual} with blocks
\begin{equation}\label{eq:ds_ns}
\begin{array}{c}
B_{0,k}=\begin{pmatrix}
a_{k}+c_{n,k} & 0 & 0 & 0 \smallskip\\
a_{k}+c_{4,k} & b_{k}+c_{4,k} & 0 & 0 \smallskip\\
a_{k}+c_{4,k} & b_{k}+c_{4,k} &c_{4,k} & b_{k}+c_{4,k} \smallskip\\
a_{k}+c_{4,k} & 0 & 0 & b_{k}+c_{4,k}
\end{pmatrix}, \,
B_{1,k}=\begin{pmatrix}
b_{k}+c_{n,k} & 0 & 0 & 0 \smallskip\\
 0 & 0 & 0 & 0 \smallskip\\
0 & 0 & 0 & 0 \smallskip\\
b_{k}+c_{4,k} & c_{4,k} & 0 & 0
\end{pmatrix},
\medskip \\
B_{i,k}=\begin{pmatrix}
c_{n,k} & 0 & 0 & 0 \smallskip\\
0 & 0 & 0 & 0 \smallskip\\
0 & 0 & 0 & 0 \smallskip\\
0 & 0 & 0 & 0
\end{pmatrix}, \, i=2, \ldots, n-2, \; \;
B_{n-1,k}=\begin{pmatrix}
b_{k}+c_{n,k} & 0 & 0 & 0 \smallskip\\
b_{k}+c_{4,k} & 0 & 0 & c_{4,k} \smallskip\\
0 & 0 & 0 & 0 \smallskip\\
0 & 0 & 0 & 0
\end{pmatrix}.
\end{array}
\end{equation}
Using \eqref{eq:mask_ds} and \eqref{eq:mask_ds_ns}, we verify that the stationary and non-stationary subdivision schemes are asymptotically equivalent of order 1. To see it, we use the Lagrange form of the remainder of the Taylor expansion to write
$$
\cos(2^{-k}h)=1-\frac{h^2}{2} 2^{-2k}+\frac{h^4}{24} 2^{-4k} \cos(\xi), \quad \xi \in (0,2^{-k}h),
$$
and
$$
\cos^2(2^{-k}h)=1-h^2 2^{-2k} +\frac{h^4}{3} 2^{-4k} \cos(2\tilde{\xi}), \quad \tilde{\xi} \in (0,2^{-k}h).
$$
The previous expression allows us to get the bounds
$$
\begin{array}{l}
|a_k-\frac{1}{2}| \leq \frac{\mathpzc{A}}{4^k}, \quad
|b_k-\frac{1}{8}| \leq \frac{\mathpzc{B}}{4^k}, \quad
|c_{n,k}-\frac{1}{4n}| \leq \frac{n^{-1} \mathpzc{C}}{4^k}, \quad \forall n \geq 4,
\end{array}
$$
with $\mathpzc{A},\mathpzc{B},\mathpzc{C}$ finite positive constants independent of $n$ and $k$.\\
The latter bounds can then be used to show that
$$
\begin{array}{lll}
\|{\cal S}_{\bc^{(k)}}-{\cal S}_{\bc}\|_{\infty}&=&\left|a_k+c_{4,k}-\frac{9}{16}\right|+2\left|b_k+c_{4,k}-\frac{3}{16}\right|+\left|c_{4,k}-\frac{1}{16}\right| \smallskip \\
&\leq& \left|a_k-\frac{1}{2}\right|+2\left|b_k-\frac{1}{8}\right|+4\left|c_{4,k}-\frac{1}{16}\right|\smallskip\\
&\leq& \frac{\mathpzc{A}+2\mathpzc{B}+\mathpzc{C}}{4^k},
\end{array}
$$
and therefore prove the asymptotical equivalence of order $1$. Indeed, using \eqref{normSO}, we arrive at
$$
\sum_{k=1}^{+\infty} 2^k \|{\cal S}_{\bc^{(k)}}-{\cal S}_{\bc}\|_{\infty} \leq (\mathpzc{A}+2\mathpzc{B}+\mathpzc{C}) \, \sum_{k=1}^{+\infty} \frac{1}{2^k}<+\infty.
$$
Summarizing, assumptions $({\it i})-({\it iii})$ of Theorem \ref{teo:G1} and assumption ({\it ii}) of Theorem \ref{theo:convergenceNEW} are satisfied.
Next, we show that $\|S_k-S\|_{\infty}\leq \frac{\mathpzc{M}}{4^k}$ for all $k \geq 1$, $n\geq 5$ and $h \in \left[0,\frac{\pi}{3}\right)$, with $\mathpzc{M}$ a finite positive constant. Indeed, by  \eqref{eq:ds} and \eqref{eq:ds_ns} we have
$$
\|S_k-S\|_{\infty} \leq  \|B_{0,k}-B_0\|_{\infty} + \|B_{1,k}-B_1\|_{\infty} +\sum_{i=2}^{n-2}\|B_{i,k}-B_i\|_{\infty} + \|B_{n-1,k}-B_{n-1}\|_{\infty}.
$$
Since
{\small
$$
\begin{array}{lll}
\|B_{0,k}-B_0\|_{\infty} &=& \max \Big \{|a_k+c_{n,k}-(\frac{1}{4n}+\frac{1}{2})|, \, |a_k+c_{4,k}-\frac{9}{16}|+|b_k+c_{4,k}-\frac{3}{16}|, \,\\
&& |a_k+c_{4,k}-\frac{9}{16}|+2|b_k+c_{4,k}-\frac{3}{16}|+|c_{4,k}-\frac{1}{16}| \Big \}\\
&=& \max \Big \{|a_k+c_{n,k}-(\frac{1}{4n}+\frac{1}{2})|, \\  && |a_k+c_{4,k}-\frac{9}{16}|+2|b_k+c_{4,k}-\frac{3}{16}|+|c_{4,k}-\frac{1}{16}| \Big \},\\
&\leq& \max \Big\{ |a_k-\frac{1}{2}|+|c_{n,k}-\frac{1}{4n}|, \,
|a_k-\frac{1}{2}|+2|b_k-\frac{1}{8}|+4|c_{4,k}-\frac{1}{16}| \Big\}\\
&\leq& \frac{1}{4^k} \max \Big\{ \mathpzc{A}+n^{-1} \mathpzc{C}, \, \mathpzc{A}+2\mathpzc{B}+\mathpzc{C} \Big\}=:\frac{\mathpzc{M}_0}{4^k},\\
\|B_{1,k}-B_1\|_{\infty} &=& \|B_{n-1,k}-B_{n-1}\|_{\infty}\\
&=&\max \Big \{ |b_k+c_{n,k}-(\frac{1}{4n}+\frac{1}{8})|, \, |b_k+c_{4,k}-\frac{3}{16}|+|c_{4,k}-\frac{1}{16}|  \Big \}\\
&\leq& \max \Big \{ |b_k-\frac{1}{8}|+|c_{n,k}-\frac{1}{4n}|, \,
 |b_k-\frac{1}{8}|+2|c_{4,k}-\frac{1}{16}|\Big \}\\
&\leq& \frac{1}{4^k} \max \Big \{ \mathpzc{B}+n^{-1} \mathpzc{C}, \, \mathpzc{B}+ \frac{1}{2}\mathpzc{C} \Big \}=:\frac{\mathpzc{M}_1}{4^k},\\
\|B_{i,k}-B_i\|_{\infty} &=&|c_{n,k}-\frac{1}{4n}| \leq \frac{n^{-1}\mathpzc{C}}{4^k}, \qquad i=2,...,n-2,
\end{array}
$$
}
for $n \geq 5$, we finally obtain the bound
$$
\|S_k-S\|_{\infty} \leq \frac{\mathpzc{M}_0+\mathpzc{M}_1+(1-\frac{3}{n}) \mathpzc{C}}{4^k} \leq \frac{\mathpzc{M}}{4^k},
$$
with $\mathpzc{M}:=\mathpzc{M}_0+\mathpzc{M}_1+\mathpzc{C}$ a finite positive constant independent of $n$ and $k$. In other words assumption  ({\it iii}) of Theorem \ref{theo:convergenceNEW} and assumption  ({\it iv}) of Theorem \ref{teo:G1} are satisfied. Since $S$ has a
dominant single eigenvalue $\lambda_0=1$ and a subdominant eigenvalue $0.5<\lambda_1<1$ with algebraic and geometric multiplicity $2$ (i.e., it is a double non-defective eigenvalue), all the assumptions of Theorem \ref{theo:convergenceNEW} and Theorem \ref{teo:G1} are verified with $\sigma=4$. Hence, the non-stationary version of Doo-Sabin scheme is convergent at extraordinary faces and  the limit surfaces obtained by such a scheme are normal continuous.\\
Figure \ref{fig:ds_examples} shows two application examples of the normalized version of such a scheme, where the normalization factor is introduced to obtain refined meshes that lie in the convex hull of the initial control points.

\subsection{Generalized exponential spline surfaces of order $\geq 3$}
In \cite{S01} a generalization of order-$d$ polynomial spline surfaces to quadrilateral meshes of arbitrary topology has been proposed.
For $d=4$, the refinement rules of the corresponding scheme are the rules of Catmull-Clark subdivision scheme \cite{CC78} which, in the regular regions of the mesh, can be given in terms of the subdivision mask
\begin{equation}\label{eq:mask_cc}
\bc=\begin{pmatrix}
\frac{1}{64} & \frac{1}{16} & \frac{3}{32} & \frac{1}{16} & \frac{1}{64} \smallskip \\
    \frac{1}{16}& \frac{1}{4}& \frac{3}{8}&  \frac{1}{4} &\frac{1}{16} \smallskip\\
    \frac{3}{32}& \frac{3}{8}& \frac{9}{16}& \frac{3}{8} & \frac{3}{32} \smallskip\\
    \frac{1}{16}& \frac{1}{4}& \frac{3}{8} &\frac{1}{4} &\frac{1}{16} \smallskip\\
    \frac{1}{64}& \frac{1}{16}& \frac{3}{32}& \frac{1}{16}& \frac{1}{64}
\end{pmatrix}.
\end{equation}
Differently, in the neighborhood of an extraordinary vertex of valence $n \geq 5$,  the subdivision matrix $S_k$ of the order-4 scheme is as in  \eqref{eq:Sk_block} with $\tilde{\alpha}=1-\frac{7}{4n},$ $\tilde{\bbeta}=\left(\frac{3}{2n^2}, \frac{1}{4n^2}, 0, 0, 0, 0\right)^T$, $\tilde{\bgamma}=(\frac{3}{8},\frac14,\frac{3}{32},\frac{1}{16},\frac{1}{64},\frac{1}{16})^T$ and $6 \times 6$ blocks
$$
\begin{array}{c}
\tilde{B}_{0}=\begin{pmatrix}
\frac{3}{8}&\frac{1}{16}&0&0&0&0 \smallskip\\
    \frac{1}{4}&\frac{1}{4}&0&0&0&0 \smallskip\\
    \frac{9}{16}&\frac{3}{32}&\frac{3}{32}&\frac{1}{64}&0&0 \smallskip\\
    \frac{3}{8}&\frac{3}{8}&\frac{1}{16}&\frac{1}{16}&0&0 \smallskip\\
    \frac{3}{32}&\frac{9}{16}&\frac{1}{64}&\frac{3}{32}&\frac{1}{64}&\frac{3}{32} \smallskip\\
    \frac{1}{16}&\frac{3}{8}&0&0&0&\frac{1}{16}
\end{pmatrix}, \qquad
\tilde{B}_{1}= \begin{pmatrix}
\frac{1}{16}&0&0&0&0&0 \smallskip\\
    \frac{1}{4}&0&0&0&0&0 \smallskip\\
    \frac{1}{64}&0&0&0&0&0 \smallskip\\
    \frac{1}{16}&0&0&0&0&0 \smallskip\\
    \frac{3}{32}&0&\frac{1}{64}&0&0&0 \smallskip\\
    \frac{3}{8}&0&\frac{1}{16}&0&0&0
\end{pmatrix}, \smallskip\\
\tilde{B}_i=\mathbf{0}_{6 \times 6}, \, i=2, \ldots, n-2, \qquad
\tilde{B}_{n-1}=\begin{pmatrix}
\frac{1}{16}&\frac{1}{16}&0&0&0&0 \smallskip\\
    0&0&0&0&0&0 \smallskip\\
    \frac{1}{64}&\frac{3}{32}&0&0&0&\frac{1}{64} \smallskip\\
    0&0&0&0&0&0 \smallskip\\
		0&0&0&0&0&0 \smallskip\\
		0&0&0&0&0&0
\end{pmatrix}.\end{array}
$$
\noindent
It is a well-known fact that Catmull-Clark scheme is convergent both in regular regions and in irregular regions, and the limit surface is $C^2$-continuous in the regular regions of the mesh and $C^1$ at the limit points of extraordinary vertices.
The subdivision symbol associated to the scheme is
$$c(z_1,z_2)=\frac{(z_1 + 1)^4(z_2 + 1)^4}{64},$$
which contains the factor $(1+z_1)(1+z_2)$. Thus it verifies assumption ({\it i}) of Theorem \ref{theo:convergenceNEW} and assumption ({\it i}) of Theorem \ref{teo:G1}.

\begin{figure}[h!]
\centering
\includegraphics[trim={0.25cm -0.2cm 0.25cm 0.25cm}, clip, height=0.20\textwidth]{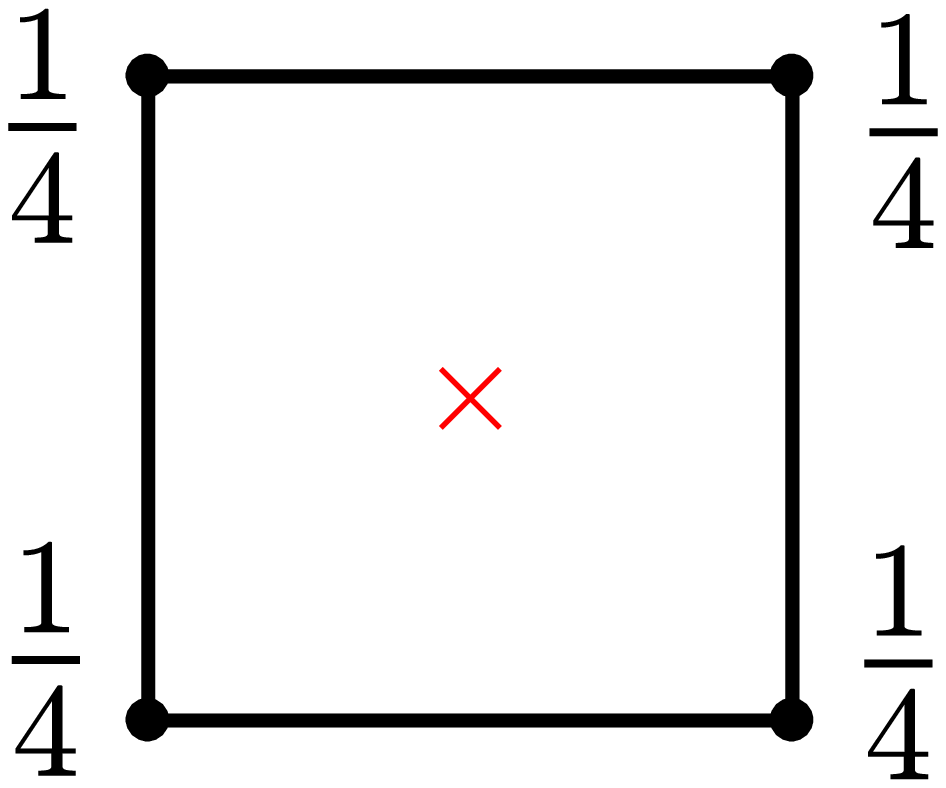}\hspace{-0.3cm}
\includegraphics[trim={0.25cm 0.25cm 0.25cm 0.25cm}, clip, height=0.33\textwidth]{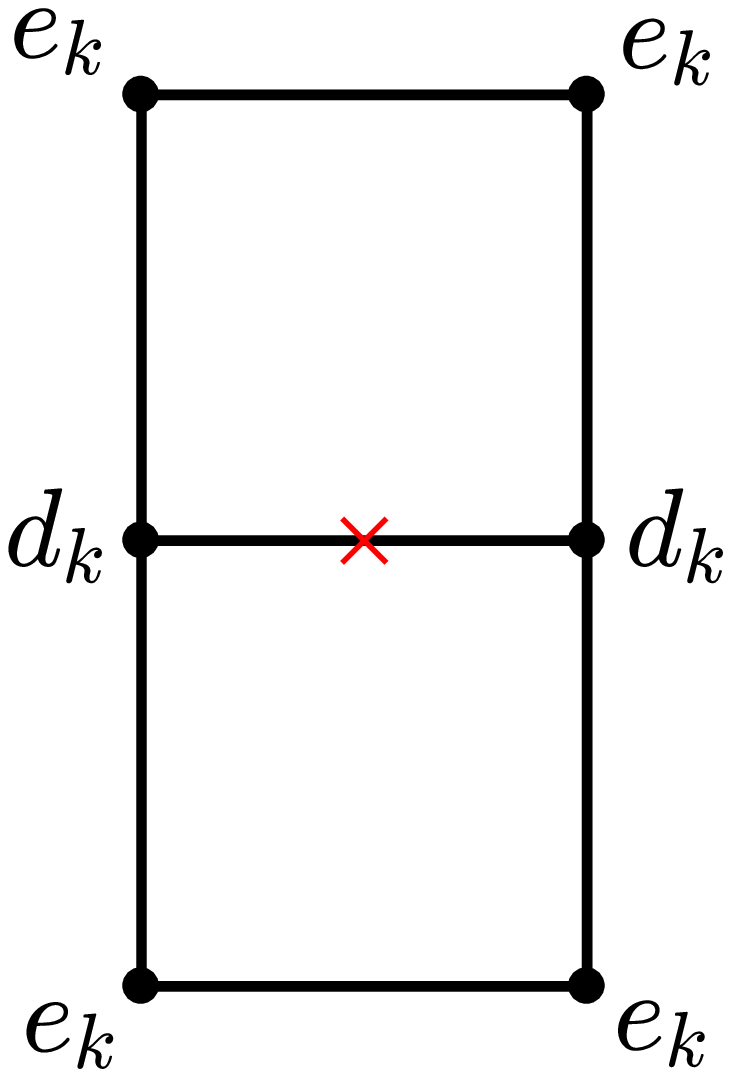}\\
\includegraphics[trim={0.25cm 0.25cm 0.25cm 0.25cm}, clip, height=0.33\textwidth]{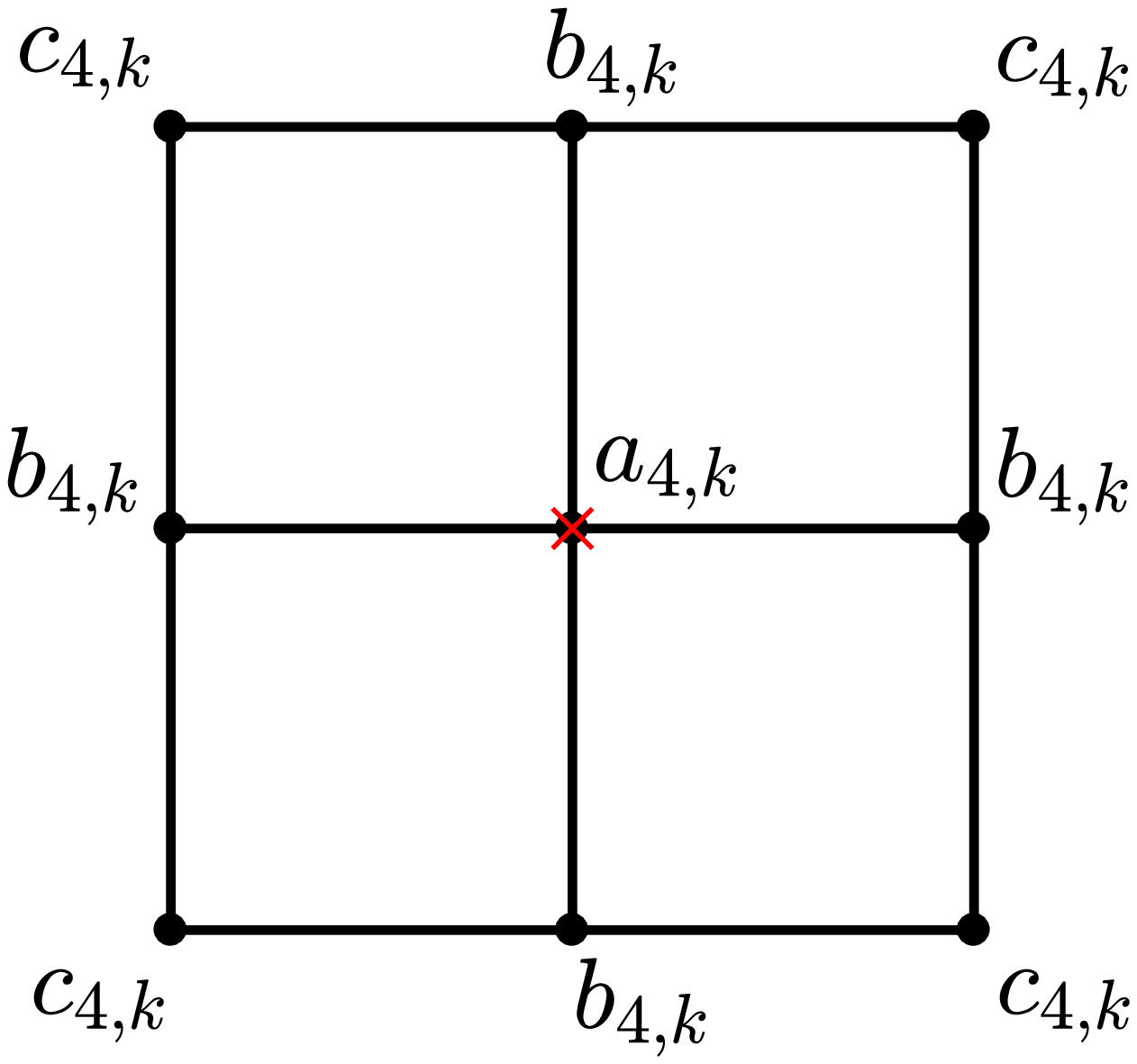}\hspace{-0.35cm}
\includegraphics[trim={0.25cm 0.25cm 0.25cm 0.15cm}, clip, height=0.38\textwidth]{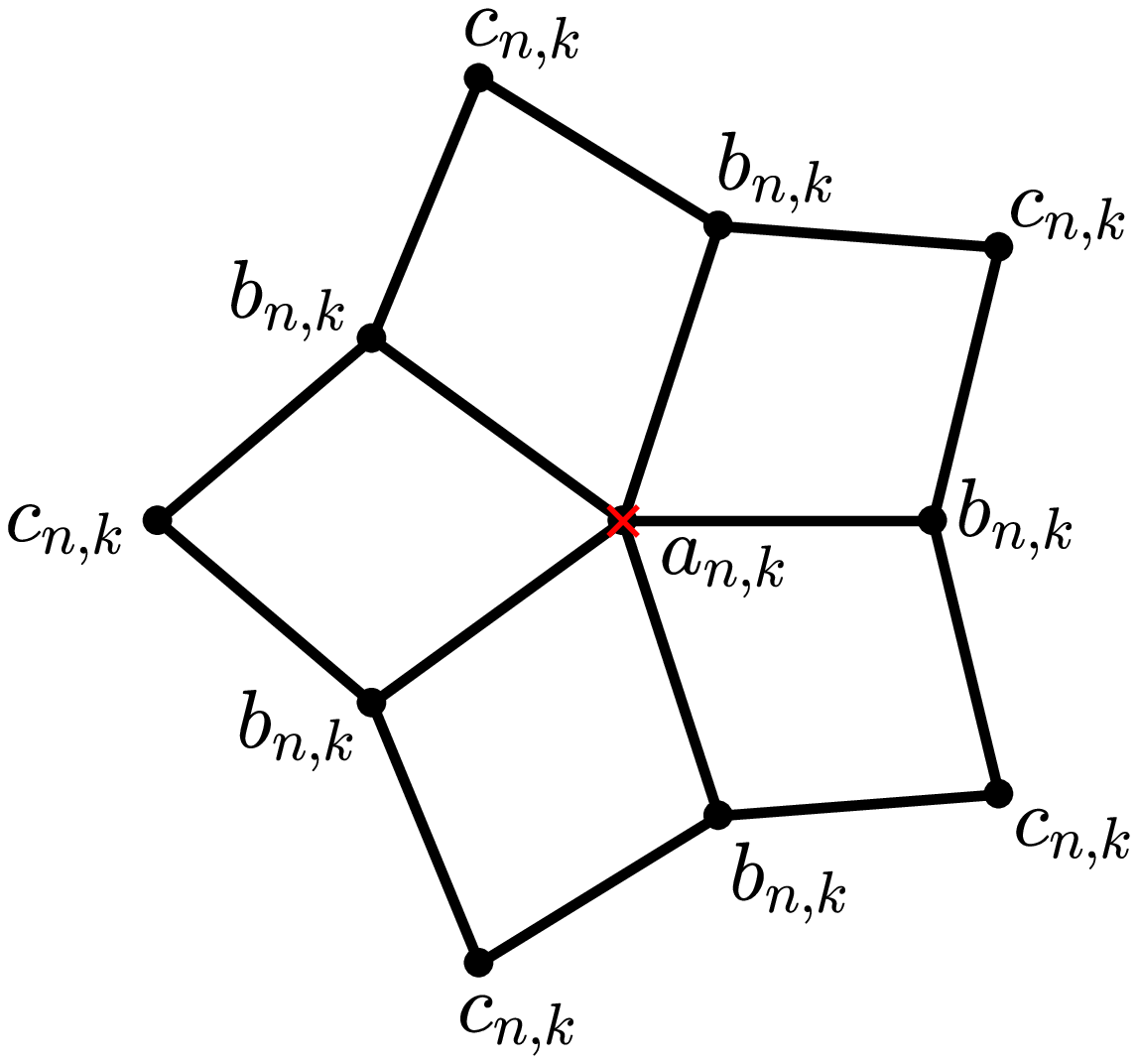}
\caption{Graphical illustration of the $k$-th level geometric refinement rules of the non-stationary subdivision scheme generalizing exponential spline surfaces of order $d=4$. Top: the red cross represents the inserted face point and edge point. Bottom: the red cross represents the inserted vertex point in the case of a regular and an extraordinary vertex. The weights appearing in the refinement rules are the ones specified in  \eqref{eq:abc_4de} and \eqref{eq:abc_n}. (Color figure online.)}
\label{fig:cc_stencils}
\end{figure}

\smallskip \noindent
The family of approximating subdivision schemes discussed in \cite{FMW14} is a non-stationary extension of the family in \cite{S01},
and provides a generalization of order-$d$ exponential spline surfaces to quadrilateral meshes of arbitrary topology.
Figure \ref{fig:cc_stencils} illustrates the $k$-th level geometric refinement rules of the order-$4$ member of this non-stationary family. We do not include a figure illustrating the topologic refinement rules since they are exactly the same as the ones used by the standard (stationary) Catmull-Clark scheme.

\begin{figure}[h!]
\centering
\subfigure[]{\includegraphics[trim={1.65cm 1.5cm 1.65cm 1.3cm}, clip, height=0.23\textwidth]{A_k1.eps}}
\subfigure[$\theta=3$]{\includegraphics[trim={2.1cm 1.8cm 2.1cm 1.4cm}, clip, height=0.23\textwidth]{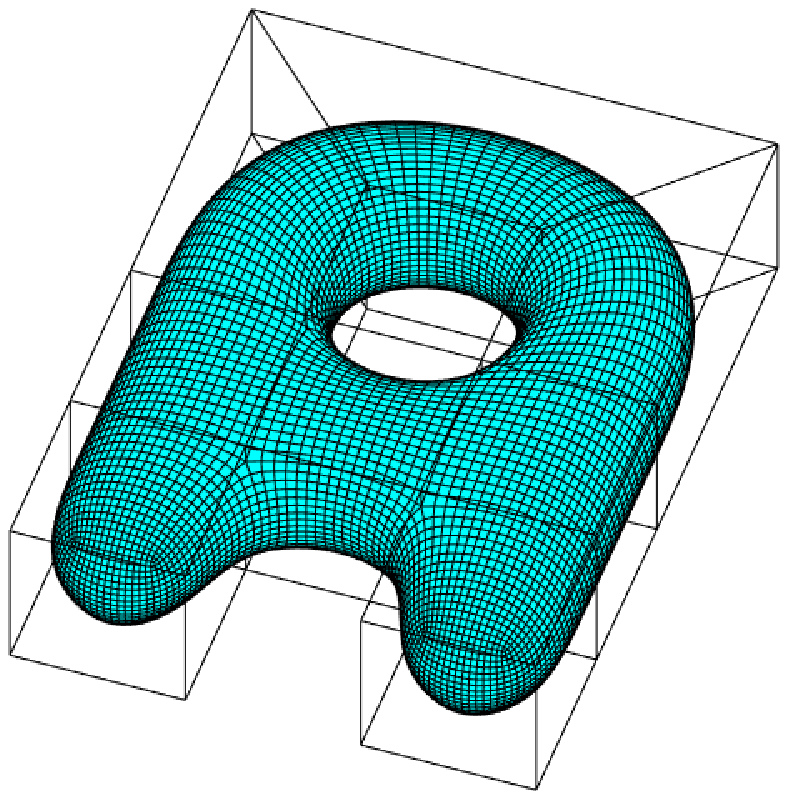}}
\subfigure[$\theta=10{\rm i}$]{\includegraphics[trim={2.1cm 1.8cm 2.1cm 1.4cm}, clip, height=0.23\textwidth]{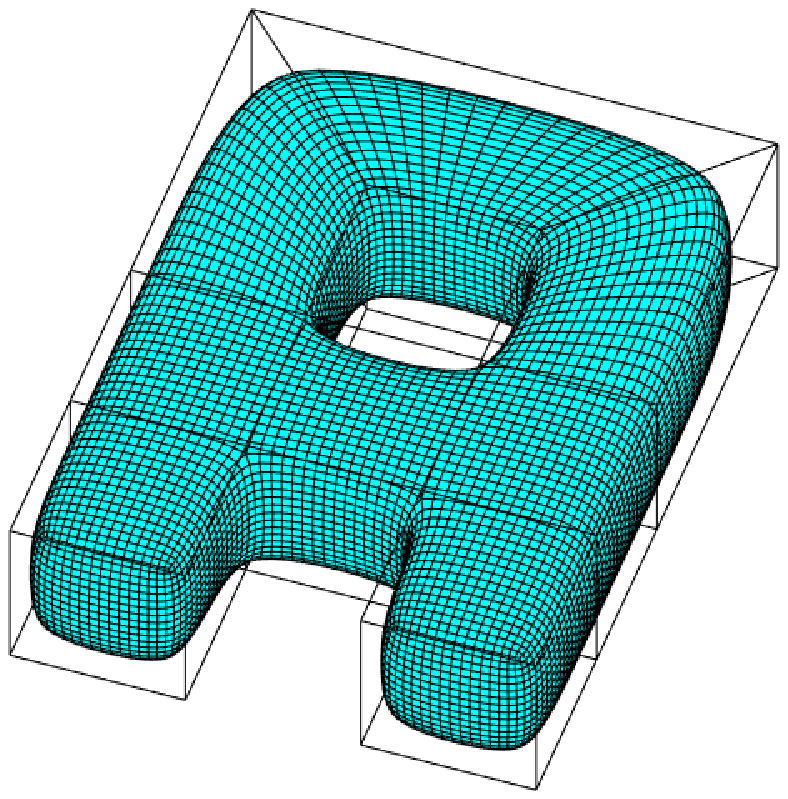}}
\vspace{0.3cm}
\setcounter{subfigure}{0}
\subfigure[]{\includegraphics[trim={2.4cm 2.1cm 2.4cm 2.1cm}, clip, height=0.22\textwidth]{triple_torus_ds_k1.eps}}
\subfigure[$\theta=3$]{\includegraphics[trim={2.8cm 2.5cm 2.8cm 2.5cm}, clip, height=0.22\textwidth]{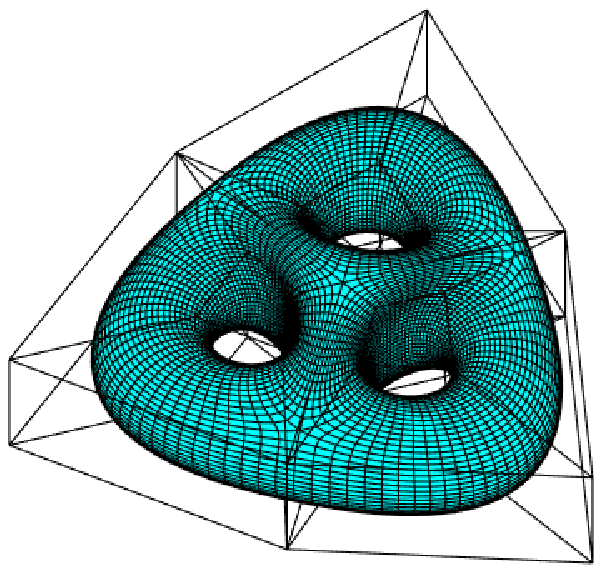}}
\subfigure[$\theta=10 {\rm i}$]{\includegraphics[trim={2.8cm 2.5cm 2.8cm 2.5cm}, clip, height=0.22\textwidth]{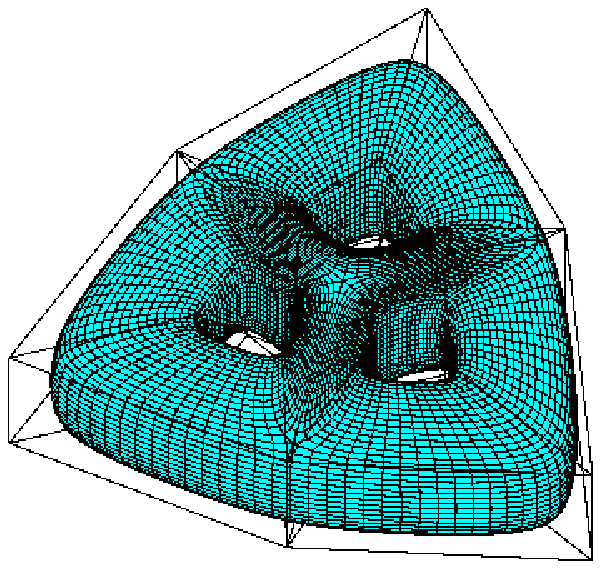}}
\caption{Original control mesh (a) and refined meshes (b,c) obtained by applying four iterations of the non-stationary subdivision scheme generalizing exponential spline surfaces of order $d=4$ with two different values of $\theta \in [0,\pi) \cup \ri (0,2 {\rm acosh}(500))$. (Color figure online.)}
\label{fig:cc_examples}
\end{figure}

\noindent
The refinement rules defining this order-$4$ non-stationary scheme, and illustrated in Figure \ref{fig:cc_stencils}, are chosen in such a way that it reproduces particular shapes such as spheres, tori or conical shapes when the initial meshes are suitably selected. In addition,  when the initial mesh is regular, the limit surface is a  tensor-product exponential spline (namely it can be either a tensor-product polynomial spline or a tensor-product trigonometric and hyperbolic spline) \cite{LWY02}.
More precisely, the $k$-th level ($k \geq 1$) refinement rules characterising the subdivision scheme depend on a $k$-th level parameter $v_{k}$ defined as
$$
v_k=\frac12 \left(e^{\frac{{\rm i}\theta}{2^{k}}}+e^{-\frac{{\rm i}\theta}{2^{k}}}\right), \quad \theta \in [0,\pi) \cup \ri (0,2 {\rm acosh}(500)), \quad k\geq 1,
$$
which satisfies
$$
(a)\; \displaystyle v_{k+1}=\sqrt{\frac{v_{k}+1}{2}},\qquad \quad (b)\; \displaystyle \lim_{k\rightarrow + \infty} v_{k}=1.
$$
Note that
$$
v_1=\frac12 \left(e^{\frac{{\rm i}\theta}{2}}+e^{-\frac{{\rm i}\theta}{2}}\right)=\left \{
\begin{array}{llll}
& \cos\left(\frac{\theta}{2}\right)\in (0,1] & {\rm if} &  \theta \in [0,\pi),\\
& \cosh\left(\frac{Im(\theta)}{2}\right)\in (1,500) & {\rm if} &  \theta \in \ri (0, 2 {\rm acosh}(500)).
\end{array}
\right.
$$
For the non-stationary approximating scheme of order $d=4$ (non-stationary version of Catmull-Clark scheme), the $k$-th level subdivision mask to be used in the regular regions of the mesh is
\begin{equation}\label{eq:mask_cc_ns}
\bc^{(k)}=\begin{pmatrix}
c_{4,k} & e_k & b_{4,k} & e_k & c_{4,k}\\
    e_k& \frac{1}{4}& d_k&  \frac{1}{4} &e_k\\
    b_{4,k}& d_k& a_{4,k}& d_k & b_{4,k}\\
    e_k& \frac{1}{4}& d_k &\frac{1}{4} &e_k\\
    c_{4,k}& e_k& b_{4,k}& e_k& c_{4,k}
\end{pmatrix},\end{equation}
where
\begin{equation}\label{eq:abc_4de}
\begin{array}{c}
\displaystyle a_{4,k}=\frac{(2v_k + 1)^2}{4(v_k + 1)^2} , \qquad b_{4,k}=\frac{2(2v_k + 1)}{16(v_k + 1)^2},\qquad c_{4,k}= \frac{1}{16(v_k + 1)^2}, \smallskip\\
\displaystyle d_k= \frac{2v_k+1}{4(v_k+1)},\qquad e_k= \frac{1}{8(v_k + 1)}.
\end{array}
\end{equation}
Hence, the associated symbol reads as
$$
c^{(k)}(z_1,z_2)=\frac{(z_1 + 1)^2(z_2 + 1)^2(z_1e^{{\rm i}\frac{\theta}{2^{k}}} + 1)(z_1 + e^{{\rm i}\frac{\theta}{2^{k}}})(z_2e^{{\rm i}\frac{\theta}{2^{k}}} + 1)(z_2 + e^{{\rm i}\frac{\theta}{2^{k}}})}{4(e^{{\rm i}\frac{\theta}{2^{k}}} + 1)^4},
$$
and contains the factor $(1+z_1)(1+z_2)$.\\
Differently, the $k$-th level subdivision matrix $\tilde{S}_k$, defined near an extraordinary vertex of valence $n\geq 5$, is of the form  \eqref{eq:Sk_block} with
\begin{equation}\label{eq:abc_n}
\begin{array}{c}
\tilde{\alpha}_k= a_{n,k},\quad
\tilde{\bbeta}_k= \displaystyle\left(b_{n,k}, c_{n,k}, 0, 0, 0, 0\right)^T, \quad
\tilde{\bgamma}_k= \displaystyle \left(d_k,\frac14,b_{4,k},e_k,c_{4,k},e_k\right)^T,\\
\displaystyle a_{n,k}=1-n(b_{n,k}+c_{n,k}), \quad b_{n,k}=\frac{2(2v_k+1)}{n^2(v_k+1)^2}, \quad  c_{n,k}=\frac{1}{n^2(v_k+1)^2},
\end{array}
\end{equation}
and $6 \times 6$ blocks
$$
\begin{array}{c}
\tilde{B}_{0,k}=\begin{pmatrix}
d_k&e_k&0&0&0&0\\
    \frac{1}{4}&\frac{1}{4}&0&0&0&0\\
    a_{4,k}&b_{4,k}&b_{4,k}&c_{4,k}&0&0\\
    d_k&d_k&e_k&e_k&0&0\\
    b_{4,k}&a_{4,k}&c_{4,k}&b_{4,k}&c_{4,k}&b_{4,k}\\
    e_k&d_k&0&0&0&e_k
\end{pmatrix},
\qquad \tilde{B}_{1,k}= \begin{pmatrix}
e_k&0&0&0&0&0\\
    \frac{1}{4}&0&0&0&0&0\\
    c_{4,k}&0&0&0&0&0\\
    e_k&0&0&0&0&0\\
    b_{4,k}&0&c_{4,k}&0&0&0\\
    d_k&0&e_k&0&0&0
\end{pmatrix}, \smallskip\\
\tilde{B}_{i,k}=\mathbf{0}_{6 \times 6}, \; i=2, \ldots, n-2, \qquad
\tilde{B}_{n-1,k}=\begin{pmatrix}
e_k&e_k&0&0&0&0\\
    0&0&0&0&0&0\\
    c_{4,k}&b_{4,k}&0&0&0&c_{4,k}\\
    0&0&0&0&0&0\\
		0&0&0&0&0&0\\
		0&0&0&0&0&0
\end{pmatrix}.
\end{array}
$$
The choice of $v_{k}$ specifies the kind of spline surface we get in the limit, in the regular regions of the mesh. In fact, if $v_{k}<1$ the scheme yields trigonometric splines, if $v_{k}=1$  polynomial splines and if $v_{k}>1$  hyperbolic splines.\\
In \cite{FMW14}, the authors prove that the limit surface obtained by applying the generalized spline schemes of order $d$ to a regular mesh is $C^{d-2}$-continuous, while in the neighborhood of extraordinary elements the $C^1$-continuity of the limit surface is shown only by numerical evidence. Here we use Theorem  \ref{theo:convergenceNEW} and Theorem \ref{teo:G1} to prove convergence and normal continuity of the limit surfaces. \\
To prove that the non-stationary version of Catmull-Clark scheme is convergent and produces normal continuous surfaces at the limit points of extraordinary vertices, we first show that  the subdivision masks $\bc$ and $\bc_k$ in \eqref{eq:mask_cc} and \eqref{eq:mask_cc_ns} are asymptotically equivalent of order 1.
To this purpose we again write
$$
\begin{array}{c}
\cos(2^{-k}\theta)=1-\frac{\theta^2}{2} 2^{-2k}+\frac{\theta^4}{24} 2^{-4k} \cos(\xi), \quad \xi \in (0,2^{-k}\theta),\\
\cos^2(2^{-k}\theta)=1-\theta^2 2^{-2k} +\frac{\theta^4}{3} 2^{-4k} \cos(2\tilde{\xi}), \quad \tilde{\xi} \in (0,2^{-k}\theta),\\
\end{array}
$$
and
$$
\begin{array}{c}
\cosh(2^{-k}\theta)=1+ \theta^2 2^{-2k} + \frac{\theta^4}{24} 2^{-4k} \cosh(\eta), \quad \eta \in (0,2^{-k}\theta),\\
\cosh^2(2^{-k}\theta)=1+ 2 \theta^2 2^{-2k} + \frac{\theta^4}{3} 2^{-4k} \cosh(2\tilde{\eta}), \quad \tilde{\eta} \in (0,2^{-k}\theta),
\end{array}
$$
from which we obtain
$$
\begin{array}{l}
|a_{4,k}-\frac{9}{16}| \leq \frac{\mathpzc{A}}{4^k}, \;
|b_{4,k}-\frac{3}{32}| \leq \frac{\mathpzc{B}}{4^k}, \;
|c_{4,k}-\frac{1}{64}| \leq \frac{\mathpzc{C}}{4^k}, \;
|d_k-\frac{3}{8}| \leq \frac{\mathpzc{D}}{4^k}, \;
|e_k-\frac{1}{16}| \leq \frac{\mathpzc{E}}{4^k}
\end{array}
$$
with $\mathpzc{A},\mathpzc{B},\mathpzc{C},\mathpzc{D},\mathpzc{E}$ finite positive constants independent of $n$ and $k$.
Thus, we get
$$
\begin{array}{lll}
\|{\cal S}_{\bc^{(k)}}-{\cal S}_{\bc}\|_{\infty}&=& \max \Big \{ |a_{4,k}-\frac{9}{16}|+4|b_{4,k}-\frac{3}{32}|+4|c_{4,k}-\frac{1}{64}|,\\
&& \qquad \qquad 2|d_k-\frac{3}{8}|+4 |e_k-\frac{1}{16}|  \Big \} \smallskip \\
&\leq& \frac{1}{4^k} \max \left \{ \mathpzc{A}+4\mathpzc{B}+4\mathpzc{C}, \, 2\mathpzc{D}+4\mathpzc{E} \right \},
\end{array}
$$
so that
$$
\sum_{k=1}^{+\infty} 2^k \|{\cal S}_{\bc^{(k)}}-{\cal S}_{\bc}\|_{\infty} \leq \max \left \{ \mathpzc{A}+4\mathpzc{B}+4\mathpzc{C}, \, 2\mathpzc{D}+4\mathpzc{E} \right \} \, \sum_{k=1}^{+\infty} \frac{1}{2^k}<+\infty.
$$
As a consequence, assumptions ({\it i})-({\it iii}) of Theorem \ref{teo:G1} and  assumption ({\it ii}) of Theorem \ref{theo:convergenceNEW} are satisfied.\\
Next, we use  formula \eqref{eq:matrix_S} to transform the matrices $\tilde{S}$ and $\tilde{S}_k$ in the block-circulant matrices denoted by ${S}$ and ${S}_k$, and verify the existence  of a finite positive constant $\mathpzc{M}$  independent of $n$ and $k$ such that  $\|S_k-S\|_{\infty}\leq \frac{\mathpzc{M}}{4^k}$ for all $k \geq 1$, $n\geq 5$ and $\theta \in [0,\pi) \cup \ri (0,2 {\rm acosh}(500))$. As before, we first write
$$
\|S_k-S\|_{\infty} \leq  \|{B}_{0,k}-{B}_0\|_{\infty} + \|{B}_{1,k}-{B}_1\|_{\infty} +\sum_{i=2}^{n-2}\|{B}_{i,k}-{B}_i\|_{\infty} + \|{B}_{n-1,k}-{B}_{n-1}\|_{\infty},
$$
and explicitly compute each norm on the right hand side as
{\small
$$
\begin{array}{lll}
\|{B}_{0,k}-{B}_0\|_{\infty} &=& \left \| \begin{pmatrix} \frac{\tilde{\alpha}_k-\tilde{\alpha}}{n} & \tilde{\bbeta}_k^T-\tilde{\bbeta}^T \smallskip\\ \frac{\tilde{\bgamma}_k-\tilde{\bgamma}}{n} & \tilde{B}_{0,k}-\tilde{B}_{0} \end{pmatrix}  \right \|_{\infty}\smallskip\\
&=& \max \{
\frac{1}{n^2} \, |\frac{7}{4}-\frac{4v_k + 3}{(v_k + 1)^2}|
+\frac{1}{n^2} \, |\frac{2(2v_k+1)}{(v_k+1)^2}- \frac{3}{2}|
+\frac{1}{n^2} \, |\frac{1}{(v_k+1)^2}-\frac{1}{4}|,\smallskip\\
&& (\frac{1}{n}+1) |d_k-\frac{3}{8}|+|e_k-\frac{1}{16}|, \smallskip\\
&& (\frac{1}{n}+2) |b_{4,k}-\frac{3}{32}|+|a_{4,k}-\frac{9}{16}|+|c_{4,k}-\frac{1}{64}|, \smallskip \\
&& (\frac{1}{n}+2) |e_k-\frac{1}{16}|+2 |d_k-\frac{3}{8}|, \smallskip \\
&& (\frac{1}{n}+2) |c_{4,k}-\frac{1}{64}|+|a_{4,k}-\frac{9}{16}|+3|b_{4,k}-\frac{3}{32}| \}, \smallskip \\
\|{B}_{1,k}-{B}_1\|_{\infty} &=& \left \| \begin{pmatrix} \frac{\tilde{\alpha}_k-\tilde{\alpha}}{n} & \tilde{\bbeta}_k^T-\tilde{\bbeta}^T \smallskip\\ \frac{\tilde{\bgamma}_k-\tilde{\bgamma}}{n} & \tilde{B}_{1,k}-\tilde{B}_{1} \end{pmatrix}  \right \|_{\infty}\smallskip \\
&=& \max \{
\frac{1}{n^2} \, |\frac{7}{4}-\frac{4v_k + 3}{(v_k + 1)^2}|
+\frac{1}{n^2} \, |\frac{2(2v_k+1)}{(v_k+1)^2}- \frac{3}{2}|
+\frac{1}{n^2} \, |\frac{1}{(v_k+1)^2}-\frac{1}{4}|,\smallskip\\
&& \frac{1}{n}|d_k-\frac{3}{8}|+|e_k-\frac{1}{16}|, \smallskip\\
&& \frac{1}{n}|b_{4,k}-\frac{3}{32}|+|c_{4,k}-\frac{1}{64}|, \smallskip\\
&& (\frac{1}{n}+1) |c_{4,k}-\frac{1}{64}|+|b_{4,k}-\frac{3}{32}|, \smallskip \\
&& (\frac{1}{n}+1) |e_k-\frac{1}{16}|+|d_k-\frac{3}{8}| \}, \smallskip\\
\|{B}_{i,k}-{B}_i\|_{\infty} &=& \left \| \begin{pmatrix} \frac{\tilde{\alpha}_k-\tilde{\alpha}}{n} & \tilde{\bbeta}_k^T-\tilde{\bbeta}^T \smallskip\\ \frac{\tilde{\bgamma}_k-\tilde{\bgamma}}{n} & \tilde{B}_{i,k}-\tilde{B}_{i} \end{pmatrix}  \right \|_{\infty}
\smallskip \\
&=& \max \{
\frac{1}{n^2} \, |\frac{7}{4}-\frac{4v_k + 3}{(v_k + 1)^2}|
+\frac{1}{n^2} \, |\frac{2(2v_k+1)}{(v_k+1)^2}- \frac{3}{2}|
+\frac{1}{n^2} \, |\frac{1}{(v_k+1)^2}-\frac{1}{4}|,\smallskip\\
&& \frac{1}{n} |d_k-\frac{3}{8}|, \, \frac{1}{n}|b_{4,k}-\frac{3}{32}|, \, \frac{1}{n} |e_k-\frac{1}{16}|, \, \frac{1}{n} |c_{4,k}-\frac{1}{64}| \}, \\
&& i=2,...,n-2, \smallskip\\
\|{B}_{n-1,k}-{B}_{n-1}\|_{\infty} &=& \left \| \begin{pmatrix} \frac{\tilde{\alpha}_k-\tilde{\alpha}}{n} & \tilde{\bbeta}_k^T-\tilde{\bbeta}^T \smallskip\\ \frac{\tilde{\bgamma}_k-\tilde{\bgamma}}{n} & \tilde{B}_{n-1,k}-\tilde{B}_{n-1} \end{pmatrix}  \right \|_{\infty}\smallskip \\
&=& \max \{
\frac{1}{n^2} \, |\frac{7}{4}-\frac{4v_k + 3}{(v_k + 1)^2}|
+\frac{1}{n^2} \, |\frac{2(2v_k+1)}{(v_k+1)^2}- \frac{3}{2}|
+\frac{1}{n^2} \, |\frac{1}{(v_k+1)^2}-\frac{1}{4}|,\smallskip\\
&& \frac{1}{n}|d_k-\frac{3}{8}|+2|e_k-\frac{1}{16}|, \smallskip \\
&& (\frac{1}{n}+1)|b_{4,k}-\frac{3}{32}|+2|c_{4,k}-\frac{1}{64}|, \frac{1}{n} |e_k-\frac{1}{16}|, \frac{1}{n} |c_{4,k}-\frac{1}{64}| \}.
\end{array}
$$
}
Moreover, in view of the bounds
$$
\begin{array}{lll}
|\frac{7}{4}-\frac{4v_k + 3}{(v_k + 1)^2}|&=&16 |\frac{7}{64}-\frac{4v_k + 3}{16(v_k + 1)^2}|  \\
&\leq& 16 \left( |\frac{3}{32}-\frac{4v_k + 2}{16(v_k + 1)^2}| + |\frac{1}{64}-\frac{1}{16(v_k + 1)^2}|  \right)\\
&\leq& \frac{16(\mathpzc{B}+\mathpzc{C})}{4^k}, \smallskip \\
|\frac{2(2v_k+1)}{(v_k+1)^2}- \frac{3}{2}|&=& 16|\frac{2(2v_k+1)}{16(v_k+1)^2}- \frac{3}{32}| \leq \frac{16 \mathpzc{B}}{4^k} ,\smallskip \\
|\frac{1}{(v_k+1)^2}-\frac{1}{4}|&=&16|\frac{1}{16(v_k+1)^2}-\frac{1}{64}| \leq \frac{16 \mathpzc{C}}{4^k},\smallskip \\
\end{array}
$$
we are finally able to bound the norms of the blocks as
{\small
$$
\begin{array}{lll}
\|{B}_{0,k}-{B}_0\|_{\infty} &\leq& \max \{
\frac{32 n^{-2}(\mathpzc{B}+\mathpzc{C})}{4^k}, \, \frac{(n^{-1}+1)\mathpzc{D}+\mathpzc{E}}{4^k}, \, \frac{(n^{-1}+2)\mathpzc{B}+\mathpzc{A}+\mathpzc{C}}{4^k}, \,
\frac{(n^{-1}+2)\mathpzc{E}+2\mathpzc{D}}{4^k} \\
&& \frac{(n^{-1}+2)\mathpzc{C}+\mathpzc{A}+3\mathpzc{B}}{4^k} \} \smallskip \\
&\leq & \frac{ \max \{ 32 (\mathpzc{B}+\mathpzc{C}), \, 2\mathpzc{D}+\mathpzc{E}, \, 3\mathpzc{B}+\mathpzc{A}+\mathpzc{C}, \, 3\mathpzc{E}+2\mathpzc{D}, \, 3\mathpzc{C}+\mathpzc{A}+3\mathpzc{B} \}}{4^k}=:\frac{\mathpzc{M}_0}{4^k}, \smallskip \\
\|{B}_{1,k}-{B}_1\|_{\infty} &\leq& \max \{
\frac{32 n^{-2}(\mathpzc{B}+\mathpzc{C})}{4^k}, \, \frac{n^{-1}\mathpzc{D}+\mathpzc{E}}{4^k}, \, \frac{n^{-1}\mathpzc{B}+\mathpzc{C}}{4^k}, \, \frac{(n^{-1}+1)\mathpzc{C}+\mathpzc{B}}{4^k}, \,
\frac{(n^{-1}+1)\mathpzc{E}+\mathpzc{D}}{4^k} \}, \smallskip\\
&\leq& \frac{ \max \{32(\mathpzc{B}+\mathpzc{C}), \, \mathpzc{D}+\mathpzc{E}, \, \mathpzc{B}+\mathpzc{C}, \, 2\mathpzc{C}+\mathpzc{B}, \, 2\mathpzc{E}+\mathpzc{D} \}}{4^k}=:\frac{\mathpzc{M}_1}{4^k}, \smallskip\\
\|{B}_{i,k}-{B}_i\|_{\infty} &\leq& \max \{
\frac{32 n^{-2}(\mathpzc{B}+\mathpzc{C})}{4^k}, \, \frac{n^{-1} \mathpzc{D}}{4^k}, \, \frac{n^{-1}\mathpzc{B}}{4^k}, \, \frac{n^{-1}\mathpzc{E}}{4^k}, \, \frac{n^{-1}\mathpzc{C}}{4^k} \}\smallskip\\
&\leq& \frac{ n^{-1} \max \{ 32(\mathpzc{B}+\mathpzc{C}), \, \mathpzc{D}, \, \mathpzc{B}, \, \mathpzc{E}, \, \mathpzc{C} \}}{4^k}=:
\frac{n^{-1} \mathpzc{M}_2}{4^k}, \quad i=2,...,n-2, \smallskip\\
\|{B}_{n-1,k}-{B}_{n-1}\|_{\infty} &\leq& \max \{
\frac{32 n^{-2}(\mathpzc{B}+\mathpzc{C})}{4^k}, \, \frac{n^{-1}\mathpzc{D}+2\mathpzc{E}}{4^k}, \, \frac{(n^{-1}+1)\mathpzc{B}+2\mathpzc{C}}{4^k}, \, \frac{n^{-1}\mathpzc{E}}{4^k},
\frac{n^{-1} \mathpzc{C}}{4^k} \} \smallskip \\
&\leq& \frac{ \max \{32 (\mathpzc{B}+\mathpzc{C}), \, \mathpzc{D}+2\mathpzc{E}, \, 2\mathpzc{B}+2\mathpzc{C}, \, \mathpzc{E}, \mathpzc{C} \} }{4^k}=:\frac{\mathpzc{M}_3}{4^k},
\end{array}
$$
}
Hence, for all $ n \geq 5$,
$$
\|S_k-S\|_{\infty} \leq  \frac{\mathpzc{M}_0+\mathpzc{M}_1+(1-\frac{3}{n})\mathpzc{M}_2+\mathpzc{M}_3}{4^k} \leq \frac{\mathpzc{M}}{4^k},
$$
with $\mathpzc{M}:=\mathpzc{M}_0+\mathpzc{M}_1+\mathpzc{M}_2+\mathpzc{M}_3$  a finite positive constant independent of $n$ and $k$.\\
The above proves that ({\it iii}) of Theorem \ref{theo:convergenceNEW} is satisfied. Moreover,  since $S$ has a
dominant single eigenvalue $\lambda_0=1$ and a subdominant eigenvalue $0.5<\lambda_1<1$ with algebraic and geometric multiplicity $2$ (double non defective eigenvalue), ({\it iv}) of Theorem \ref{teo:G1} is also satisfied with
$\sigma=4$.
It follows that all the assumptions of Theorem \ref{theo:convergenceNEW} and Theorem \ref{teo:G1} are verified. Thus, this non-stationary version of Catmull-Clark scheme is convergent at extraordinary vertices and  the limit surfaces obtained by such a scheme are normal continuous at the limit points of extraordinary vertices. Figure \ref{fig:cc_examples} shows two application examples of such a scheme.

\section*{Acknowledgments}
This research has been accomplished within RITA (Research ITalian network on Approximation).
The authors are members of the INdAM Research group GNCS, which has partially supported this work.
The authors wish to thank the reviewers for their constructive comments that allowed them to improve the presentation of the results and to clarify all the details of their proofs.

\bibliographystyle{siamplain}

\end{document}